\documentclass[graybox]{svmult}

\usepackage{type1cm}        
\usepackage{makeidx}         
\usepackage{graphicx}        
\usepackage{multicol}        
\usepackage[bottom]{footmisc}

\usepackage{newtxtext}       %
\usepackage[varvw]{newtxmath}       

\usepackage{algorithm}
\usepackage{algpseudocode}

\usepackage[shortlabels]{enumitem}
\usepackage{comment}
\usepackage{todonotes}

\newtheorem{assumption}{Assumption}

\DeclareMathOperator*{\argmin}{\arg\!\min}
\newcommand{\eig}{\tau} 

\newcommand*{\dg}[1]{\textcolor{cyan}{#1}} 
\newcommand*{\todoDG}[1]{\textcolor{red}{#1}} 

\makeindex             

\begin{document}

\title*{Approximating partial differential equations without boundary conditions}
\author{Andrea Bonito and Diane Guignard}
\institute{Andrea Bonito \at Department of Mathematics, Texas A\&M University, College Station, TX 77845, USA, \email{bonito@tamu.edu}
\and Diane Guignard \at Department of Mathematics and Statistics, University of Ottawa, ON K1N 6N5, Canada, \email{dguignar@uottawa.ca}}
%
%
\maketitle

\abstract{We consider the problem of numerically approximating the solutions to an elliptic partial differential equation (PDE) for which the boundary conditions are lacking. 
To alleviate this missing information, we assume to be given measurement functionals of the solution. 
In this context, a near optimal recovery algorithm based on the approximation of the Riesz representers of  these functionals in some intermediate Hilbert spaces is proposed and analyzed in \cite{BBCDDP2023}. Inherent to this algorithm is the computation of $H^s$, $s>1/2$, inner products on the boundary of the computational domain. We take advantage of techniques borrowed from the analysis of fractional diffusion problems to design and analyze a fully practical near optimal algorithm not relying on the challenging computation of $H^s$ inner products.
}

\section{Introduction}
\label{sec:intro}

In many applications, a lack of full information on boundary conditions arises for various reasons \cite{xu2021explore,brunton2020machine,duraisamy2019turbulence}. 
For example, the correct physics might not be fully understood \cite{raj,richards2019appropriate}, the boundary values are not accessible \cite{grinberg2008outflow} or they must be appropriately modified in numerical schemes \cite{formaggia2002numerical,richards2011appropriate}. 
In that cases, measurements can be used to alleviate the missing information; see for instance \cite{grinberg2008outflow,raissi2020hidden,BBCDDP2023}. 

We consider the model elliptic partial differential equation (PDE)
\begin{equation} \label{eqn:pb}
    -\Delta u = f \quad \text{in} \,\, \Omega,
\end{equation}
satisfied by the function $u \in H^1(\Omega)$ of interest on a bounded polygonal domain $\Omega\subset\mathbb{R}^{d}$ ($d=2,3$) and where $f\in H^1_0(\Omega)^*$, the dual of $H^1_0(\Omega)$. Because the influence of the forcing term $f$ is marginal, we assume for the purpose of this introduction that $f \equiv 0$, i.e., the function of interest $u$ is harmonic in $\Omega$
\begin{equation} \label{eqn:pb_harmonic}
    -\Delta u = 0 \quad \text{in} \,\, \Omega.
\end{equation}
Notice that \eqref{eqn:pb_harmonic} must be supplemented by suitable boundary conditions to uniquely determine $u$. In lack of such boundary conditions, we assume to have access to measurements of $u$, namely to the vector ${\boldsymbol \omega} \in \mathbb R^m$ given by
\begin{equation} \label{eqn:lambda}
    {\boldsymbol \omega} = {\boldsymbol \lambda}(u):=(\lambda^1(u),\lambda^2(u),\ldots,\lambda^m(u))\in\mathbb{R}^m,
\end{equation}
where for $j=1,...,m$, $\lambda^j$ are known functionals in $H^1(\Omega)^*$, the dual of $H^1(\Omega)$.
Still, this additional information does not sufficiently narrow down the set of functions $u$ satisfying the PDE \eqref{eqn:pb_harmonic} and the measurements \eqref{eqn:lambda}. The last ingredient is a model class assumption (prior knowledge), which restricts the possible candidates to a compact subset of $H^1(\Omega)$. The model class assumption in our context are functions $v \in H^1(\Omega)$ with traces on $\Gamma:=\partial \Omega$ satisfying $v|_{\Gamma} \in K_\Gamma^s$ for some $s > 1/2$, where 
\begin{equation}
    K^s_\Gamma:=U(H^s(\Gamma)):=\{g\in H^s(\Gamma): \,\, \|g\|_{H^s(\Gamma)}\le 1\}. 
\end{equation}
To characterize further the model class, we denote by
\begin{equation} \label{def:K}
\mathcal H^s(\Omega):=\left\{u\in H^1(\Omega): \,\,  u \text{ satisfies } \eqref{eqn:pb_harmonic} \,\, \text{and} \,\, u|_{\Gamma}  \in H^s(\Gamma)\right\} \subset H^1(\Omega),
\end{equation}
the space of harmonic functions with traces in $H^s(\Gamma)$. 
The model class is the unit ball in $\mathcal H^s(\Omega)$, i.e.,
$$
K^s:=U(\mathcal H^s(\Omega)),
$$
and we gather all the information available on the function to be recovered in
\begin{equation} \label{def:K_w}
    K^s_{ {\boldsymbol \omega}}:=\{u\in K^s: \,\, {\boldsymbol \lambda}(u)={\boldsymbol \omega}\}.
\end{equation}

The \emph{optimal recovery problem} is to identify a function in $K^s_{\boldsymbol \omega}$ that would be the best representative of all $K^s_{\boldsymbol \omega}$. More precisely, given ${\boldsymbol \omega}\in\mathbb{R}^m$ and $s>1/2$, it consists of finding $u^* \in H^1(\Omega)$ that minimizes the worst-case recovery error, namely
\begin{equation} \label{def:optimal_u_star}
    u^* \in \argmin_{v\in H^1(\Omega)}\sigma(v;K^s_{{\boldsymbol \omega}}), \quad \text{where} \quad \sigma(v;K^s_{{\boldsymbol \omega}}):= \sup_{u\in K^s_{{\boldsymbol \omega}}}\|u-v\|_{H^1(\Omega)}.
\end{equation}
It is well-known \cite{MR1985} that the solution to this problem is unique and $u^* \in H^1(\Omega)$ is the center of the smallest closed ball in $H^1(\Omega)$ containing $K^s_{{\boldsymbol \omega}}$ (the Chebyshev ball).
Moreover, we have that $\sigma(u^*;K^s_{{\boldsymbol \omega}}) = R(K^s_{\boldsymbol \omega})$ is the radius of the Chebyshev ball which shall serve as a benchmark for the performances of practical algorithms used to construct approximations of $u^*$.
Equipped with the standard $H^s(\Gamma)$ inner product,  $\mathcal H^s(\Omega)$ is an Hilbert space. This together with the fact that $K^s$ is the unit ball in $\mathcal H^s(\Omega)$, we must have that  $u^*$ is the element of $K^s_\omega$ with minimal norm \cite{MR1985}. This crucial property offers a computationally favorable expression
\begin{equation}\label{e:representation_riesz}
u^* = \sum_{i=1}^m U_i \phi^i,
\end{equation}
where $\phi^i \in \mathcal H^s(\Omega)$ is the Riesz representer of the functional $\lambda^i$ over $\mathcal H^s(\Omega) \subset H^1(\Omega)$, i.e., 
\begin{equation}\label{e:riesz}
\langle \phi^i,v \rangle_{\mathcal H^s(\Omega)} = \lambda^i(v), \qquad \forall\, v \in \mathcal H^s(\Omega).
\end{equation}
The coefficients ${\bf U} :=( U_i )_{i=1}^m \subset \mathbb R^m$ are determined by the relation
$$
G {\bf U} =  {\boldsymbol \omega},
$$
where $G=(g_{ij})_{i,j=1}^m$ with $g_{ij} := \lambda^i(\phi^j)$. 
Since $\phi^i \in \mathcal H^s(\Omega)$, it is harmonic and thus $\phi^i = E \psi^i$, where $E: H^{1/2}(\Gamma) \rightarrow H^{1}(\Omega)$ is the harmonic extension operator (see \eqref{eq:E}) and $\psi^i \in H^s(\Gamma)$ satisfies
\begin{equation}\label{e:riesz_bdy}
\langle \psi^i,v \rangle_{H^s(\Gamma)} = (\lambda^i \circ E)(v), \qquad \forall\, v \in H^s(\Gamma).
\end{equation}
The algorithm proposed in \cite{BBCDDP2023} takes advantage of the representation~\eqref{e:representation_riesz}, where $\phi^i$ and $\psi^i$ are replaced by their respective finite element approximations. 

We reiterate that the algorithm in \cite{BBCDDP2023} is fully implementable provided a finite element approximation of \eqref{e:riesz} can be constructed. While standard finite element methods can be used for integral $s$, the approximation of the PDE \eqref{e:riesz_bdy} for nonintegral $s$ is computationally challenging in presence of the nonlocal $H^s(\Gamma)$ inner product. We focus on the case $1/2 < s <1$ since the general case, modulo additional technicalities, can be obtain by recurrence as described in Section~\ref{sec:num_res}. For $1/2 < s <1$, we advocate the use of the following inner product on $H^s(\Gamma)$
\begin{equation}\label{e:inner_intro}
\langle \mathcal{L}^{s/2} v,\mathcal{L}^{s/2} w\rangle_{L^2(\Gamma)}, \qquad \forall\, v,w \in H^s(\Gamma),
\end{equation}
where $\mathcal L:=(I-\Delta_\Gamma)$ and $\mathcal{L}^{s/2}$ denotes its fractional powers, see Section~\ref{sec:Hs}.
The advantage of using \eqref{e:inner_intro} as inner product is that \eqref{e:riesz_bdy} becomes a fractional diffusion problem
\begin{equation}\label{e:FD}
(I-\Delta_\Gamma)^s \dot \psi^i = \lambda^i \circ E.
\end{equation}
Hereafter we use the dot notation to indicate quantities defined using the alternate inner product \eqref{e:FD}. 

The purpose of this work is to design a finite element method for its approximation and accounts for its effect in the analysis of the optimal recovery algorithm. The approximation of fractional diffusion problems is an active research topic. 
We refer to the reviews \cite{bonito2018numerical,lischke2020fractional} and the reference therein for details. In this work, we follow the strategy initiated in \cite{BP2015} later refined in \cite{bonito2017numerical,BLP2019}. Compared to the setting considered in the previous references, the fractional diffusion problem \eqref{e:FD} is set on a Lipschitz closed hyper-surfaces $\Gamma$ and involves a forcing term $\lambda^i \circ E$ that is not necessarily in $L^2(\Gamma)$. 
Nevertheless, the strategy proposed in \cite{BP2015} relies on the following integral representation of $\dot \psi^i$ in \eqref{e:FD}
\begin{equation}\label{e:int}
\dot \psi^i= \frac{\sin(\pi s)}{\pi}\int_{-\infty}^{\infty}e^{(1-s)y}(e^yI+\mathcal{L})^{-1}(\lambda^i\circ E){\rm d}y
\end{equation}
and consists of employing a sinc quadrature for the indefinite integral and a finite element method for the approximation of each integrand at the quadrature points.
The resulting numerical method, see Algorithm~\ref{alg:practical}, is straightforward to implement as it only require a finite element solver for standard reaction-diffusion problems. Once the approximation $\dot\psi^i_{k,h}$ is computed, the approximation $\dot\phi^i_{k,h}$ is defined as the discrete harmonic extension of $\dot\psi^i_{k,h}$; see Section~\ref{sec:tools}.

The error analysis performed in this work is inspired from \cite{BGL2024}, see also \cite{BL2022}, where fractional diffusion problems on closed smooth manifolds with a specific forcing term (white noise) are addressed. Theorem~\ref{t:main} guarantees that
$$
\| \dot \phi^i - \dot \phi^i_{k,h} \|_{H^{1}(\Omega)} \in \mathcal O(e^{-\frac{\pi^2}{k}} + h^{q}),
$$
i.e., the error in computing the Riesz representers is exponentially decreasing as the sinc quadrature parameter $k$ decreases and exhibits an algebraic decay rate $h^{q}$ in terms of the finite element resolution characterized by the mesh size $h$, where $q>0$ depends on various parameters.
This result is used in Corollary~\ref{c:main} to provide a recovery error estimate for the output of Algorithm~\ref{alg:OR}
$$
\| u - u_{k,h} \|_{H^1(\Omega)} \leq R(\dot K^s_{\boldsymbol \omega}) + \mathcal O(e^{-\frac{\pi^2}{k}} + h^{q}).
$$

In Section~\ref{sec:prelim}, we introduce the optimal recovery algorithm (Algorithm~\ref{alg:OR}) from \cite{BBCDDP2023} and justify expression~\eqref{e:FD} along with the integral representation \eqref{eq:psi_Balak} of $\dot\psi^i$.
The approximation of the latter is discussed in Section~\ref{sec:algos}. 
We propose a \emph{practical} algorithm (Algorithm~\ref{alg:practical}), which build on the \emph{theoretical} algorithm (Algorithm~\ref{alg:theoretical}), but provides an efficient way to approximate the Riesz representers when $s$ is non-integral without computing demanding $H^s(\Gamma)$ inner products.
The performance of the practical algorithm is analyzed in Section~\ref{sec:error}, which contains the main results of this work (Theorem~\ref{t:main} and Corollary~\ref{c:main}).
We illustrate numerically the performance of Algorithm~\ref{alg:OR} and analyze the effect of the parameters $s$, $m$ and $h$ on the recovery error. 
We postpone to the appendix results related to the approximation of fractional diffusion problems, which may not be readily available in the literature but whose proofs are similar to existing ones. 

From now on, $C$ denotes a generic constant independent of the discretization parameters (i.e., the mesh size $h$ and the sinc quadrature spacing $k$) and we use the notation $a\lesssim b$ to indicate that $a\le Cb$.

\section{Preliminaries}\label{sec:prelim}

\subsection{The optimal recovery algorithm}\label{s:algo}

In this section, we recall the algorithm proposed \cite{BBCDDP2023} along with its properties. 
We return to the PDE \eqref{eqn:pb} with $f \in H^1_0(\Omega)^*$ and decompose $u$ as 
\begin{equation}\label{e:split}
u=u_f + u_{\mathcal H},
\end{equation}
where $u_f \in H^1_0(\Omega)$ solves 
\begin{equation}\label{e:u_rhs}
-\Delta u_f = f \qquad \textrm{in }\Omega
\end{equation}
and $u_{\mathcal H} \in U(\mathcal H^s(\Omega))$. Since $u_f \in H^1_0(\Omega)$ is uniquely determined by \eqref{e:u_rhs}, the recovery procedure only targets $u_{\mathcal H}$.
Therefore, without loss of generality, we assume from now on that $f=u_f=0$.

\begin{algorithm}
\caption{Recovery Algorithm - simplified from  \cite{BBCDDP2023}}\label{alg:OR}
\begin{algorithmic}
\Require $\varepsilon>0$
\Ensure $\widehat u$
\State $\diamond$ Compute an approximation $\widehat u_f$ to the function $u_f \in H^1_0(\Omega)$ in \eqref{e:u_rhs} such that 
\State $\| u_f - \widehat u_f \|_{H^1(\Omega)} \leq \varepsilon$
\State $\diamond$ Set $\widehat{\boldsymbol \omega}:={\boldsymbol \omega}-\lambda(\widehat u_f)$
\State $\diamond$ For each $j\in\{1,2,\ldots,m\}$, compute an approximation $\widehat \phi^j$ to the Riesz representer $\phi^j$ in \eqref{e:riesz} such that $\| \phi^j - \widehat \phi^j \|_{H^1(\Omega)} \leq \varepsilon$
\State $\diamond$ Set $\widehat G = (\widehat g_{ij})_{ij=1}^m$ with $\widehat g_{ij} =\lambda^i(\widehat\phi^j)$
\State $\diamond$ Find the coefficients $\widehat{\bf  U} = (\widehat U_i)_{i=1}^m$ satisfying $\widehat G \widehat{\bf U} = \widehat {\boldsymbol \omega}$ \Comment{Use Moore-Penrose inverse of $\widehat G$}
\State $\diamond$  Set $\widehat u = \sum_{i=1}^m \widehat U_i \widehat \phi^i + \widehat u_f$.
\end{algorithmic}
\end{algorithm}

We have already pointed out that Algorithm~\ref{alg:OR} is fully implementable provided approximations of the Riesz representers can be computed. The purpose of this work is to design and analyze an efficient algorithm for the approximation of the Riesz representers.

Theorem~3.1 in \cite{BBCDDP2023} guarantees that there is a constant $\Lambda$ depending on the measurements $\lambda^i$ and their Riesz representers $\phi^i$, $i=1,...,m$, on $m$ and ${\boldsymbol \omega}$ as well as the constant $C_s$ in
\begin{equation}\label{e:Cs}
\| v \|_{H^1(\Omega)} \leq C_s \| v \|_{\mathcal H^s(\Omega)}, \qquad \forall\, v \in \mathcal H^s(\Omega),
\end{equation}
such that
\begin{equation}\label{e:OR}
\| u- \widehat u \|_{H^1(\Omega)} \leq R(K^s_{\boldsymbol \omega}) + \Lambda \varepsilon. 
\end{equation}
In particular, for $\varepsilon$ sufficiently small, we have
\begin{equation}\label{e:OR2}
\| u- \widehat u \|_{H^1(\Omega)} \leq C(\varepsilon) R(K^s_{\boldsymbol \omega}),
\end{equation}
where $C(\varepsilon) \to 1^+$ as $\varepsilon \to 0^+$. When the output of a recovery algorithm satisfies \eqref{e:OR2}, we say that the algorithm is near-optimal.

\subsection{The space $H^s(\Gamma)$ and the model class assumption}
\label{sec:Hs}

Let $H^0(\Gamma):=L^2(\Gamma)$ be the space of square integrable functions over $\Gamma$ with inner product and associated norm
$$
\langle v,w \rangle_{L^2(\Gamma)}:= \int_\Gamma v w, \qquad \| v \|_{L^2(\Gamma)}^2:= \langle v,v \rangle_{L^2(\Gamma)}, \qquad  \forall\, v,w \in L^2(\Gamma).
$$
We shall frequently use without mentioning the identification of $L^2(\Gamma)$ with its dual.  Also, define
$$
H^1(\Gamma):= \{ v \in L^2(\Gamma): \,\, | \partial_{i} v |_{L^2(\Gamma)} < \infty, \ \forall\, i=1,...,d \},
$$
where $\partial_{i}$ are the components of the tangential gradient $\nabla_\Gamma$. It is an Hilbert space when equipped with the inner product
$$
\langle v,w \rangle_{H^1(\Gamma)}:=\langle v,w \rangle_{L^2(\Gamma)} + \sum_{i=1}^d\langle \partial_{i} v, \partial_{i} w \rangle_{L^2(\Gamma)},  \qquad \| v \|^2_{H^1(\Gamma)}:= \langle v,v \rangle_{H^1(\Gamma)}.
$$
Its dual space is denoted $H^{-1}(\Gamma)$ and $\langle .,.\rangle$ stands for the duality pairing between $H^{-1}(\Gamma)$ and $H^1(\Gamma)$.

Sobolev fractional spaces on surfaces can be defined using an intrinsic norm, see for instance \cite{BBCDDP2023}.
However, it will be more convenient to use a spectral characterization of these fractional spaces.
In this aim, we define $\mathcal{T}: H^{-1}(\Gamma)\rightarrow H^{1}(\Gamma)$, where for $g\in H^{-1}(\Gamma)$, $w:= \mathcal Tg$ is the unique function in $H^1(\Gamma)$ satisfying
\begin{equation} \label{eq:def_op_L}
\int_{\Gamma}wv+\int_{\Gamma}\nabla_{\Gamma}w\cdot\nabla_{\Gamma}v = \langle g,v\rangle,
 \qquad \forall\, v\in H^1(\Gamma).
\end{equation}
The operator $\mathcal{L}:=I-\Delta_{\Gamma}: H^1(\Gamma) \rightarrow H^{-1}(\Gamma)$ is the inverse of $\mathcal T$.

Let $\{(\eig_j,e_j)\}_{j=1}^{\infty} \subset \mathbb R \times H^1(\Gamma)$ be the eigenpairs of $\mathcal{L}$, ordered and normalized so that the eigenvalues and eigenfunctions satisfy
$$1=\eig_1<\eig_2\le \eig_3\le \cdots , \qquad \text{and} \qquad \langle e_i,e_j\rangle_{L^2(\Gamma)}=\delta_{ij} \quad \forall\, i,j\ge 1.$$
For any $g\in L^2(\Gamma)$ we thus have
\begin{equation} \label{eqn:eig_expansion}
  g=\sum_{j=1}^{\infty}g_je_j, \qquad g_j:=\langle g,e_j\rangle_{L^2(\Gamma)}, 
\end{equation}
which offers an alternate characterization of $L^2(\Gamma)$.
This motivates the introduction of $\dot H^{p}(\Gamma)$, $p \ge  0$, as the space of functions in $L^2(\Gamma)$ for which
\begin{equation} \label{def:spaceHs}
\|g\|_{\dot H^p(\Gamma)}^2= \langle g,g\rangle_{\dot H^p(\Gamma)} < \infty, 
\end{equation}
where
\begin{equation}\label{def:spaceHs_innerproduct}
\langle f,g\rangle_{\dot H^p(\Gamma)} := \sum_{j=1}^{\infty}\langle f,e_j\rangle_{L^2(\Gamma)} \langle g,e_j \rangle_{L^2(\Gamma)} \eig_j^p<\infty.
\end{equation}
We identify the duality product $\langle .,.\rangle$ between $\dot H^{-p}(\Gamma)$ and
$\dot H^{p}(\Gamma)$ with the $L^2(\Gamma)$ scalar product so that
$$
\bigg\langle \sum_{j=1}^\infty c_j e_j, \sum_{j=1}^\infty d_j e_j \bigg\rangle := \sum_{j=1}^\infty c_j d_j \qquad
     \text{for} \quad \sum_{j=1}^\infty d_j e_j \in \dot H^{p}(\Gamma).
$$
This leads to the following definitions for negative indices
\begin{equation}
     \dot H^{-p}(\Gamma):= \bigg\{ \bigg\langle \sum_{j=1}^\infty c_je_j,. \bigg\rangle
     : \ \{c_j\}_{j=1}^\infty \subset \mathbb R, \ \sum_{j=1}^\infty c_j^2\eig_j^{-p} < \infty \bigg\}, \quad p > 0.
\end{equation}
We refer to \cite{bonito2017numerical} for more details.

The regularity properties of the operators $\mathcal T$ and $\mathcal L$ are critical to assess the equivalence of the spaces $\dot H^q(\Gamma)$ with the standard intermediate spaces $H^q(\Gamma)$. We make the following assumption. 

\begin{assumption}\label{a:pickup}
We assume that there is $\alpha_{\Gamma}>0$ such that the operator $\mathcal{T}$ has elliptic regularity pickup index $\alpha_{\Gamma}$. That is, $\mathcal{T}$ is an isomorphism from $H^{-1+q}(\Gamma)$ to $H^{1+q}(\Gamma)$ for any $q\in [0,\alpha_{\Gamma}]$.
\end{assumption}

When Assumption~\ref{a:pickup} holds, a similar argument that the one used in Proposition~4.1 in \cite{BP2015}, see also \cite{bonito2017numerical}, guarantees that for $q \in  [-1-\alpha_\Gamma,1+\alpha_\Gamma]$ we have 
\begin{equation}\label{e:equiv_Hs}
\dot H^q(\Gamma) \cong H^q(\Gamma)
\end{equation}
with equivalent norms. We denote by $d_s$ and $D_s$ the constants only depending on $s$ and $\Gamma$ such that
\begin{equation}\label{e:equiv_norm}
d_s \| v \|_{\dot H^s(\Gamma)} \leq \| v \|_{H^s(\Gamma)} \leq D_s \| v \|_{\dot H^s(\Gamma)}, \qquad \forall\, v \in H^s(\Gamma).
\end{equation}

The sets $\dot K^s_\Gamma$, $\dot K^s$, and $\dot K^s_{\boldsymbol \omega}$ associated with the recovery problem are as their counterparts without the dot symbol but using the norm $\|.\|_{\dot H^s(\Gamma)}$ instead of $\|.\|_{H^s(\Gamma)}$. With these notations, the model class assumption reads
\begin{assumption}\label{ass:model}
There exists $s>1/2$ such that the function $u$ to be recovered belongs to $\dot K^s_{{\boldsymbol \omega}}$. 
\end{assumption}
With this notation, the recovery estimate \eqref{e:OR} reads
\begin{equation}\label{e:OR_dot}
\| u- \hat u \|_{H^1(\Omega)} \leq R(\dot K^s_{\boldsymbol \omega}) + \dot \Lambda \varepsilon,
\end{equation}
for some constant $\dot \Lambda$  depending on the measurements $\lambda^i$ and their Riesz representers $\dot\phi^i$, $i=1,...,m$ in $\dot H^s(\Gamma)$, on $m$ and ${\boldsymbol \omega}$ as well as the constant $C_s$ in \eqref{e:Cs}. Moreover, we have
$$
c R(\dot K^s_{\boldsymbol \omega}) \leq R(K^s_{\boldsymbol \omega}) \leq C R(\dot K^s_{\boldsymbol \omega}),
$$
where the constants $c,C$ depends on the equivalence constants in \eqref{e:equiv_norm}.

\subsection{The role of fractional diffusion} \label{sec:role_frac}

In this section, we further exploit the structure of the norm $\|.\|_{\dot H^s(\Gamma)}$ to deduce a characterization using the (spectral) fractional powers of $\mathcal L$.
In turns, this justify the fractional diffusion formulation \eqref{e:riesz} for the Riesz representers. 

Recall that for $1/2 < s < 1$, the spaces $\dot H^s(\Gamma)$ and $H^s(\Gamma)$ are the same. Also, $\{ (\tau_j,e_j) \}_{j=1}^\infty$ denotes the eigenpairs of $\mathcal L$ and define for $\beta \in (-1,1)$ and $q \in \mathbb R$, $\mathcal{L}^\beta: \dot H^{q}(\Gamma) \rightarrow \dot H^{q-2\beta}(\Gamma)$ by 
\begin{equation}\label{e:spectral}
     \mathcal{L}^\beta v := \sum_{j=1}^\infty \eig_j^\beta v_j e_j, \qquad \forall\, v=\big\langle \sum_{j=1}^\infty v_j e_j,. \big\rangle \in \dot H^q(\Gamma).
\end{equation}
Thanks to the Lax-Milgram theory and the equivalence \eqref{e:equiv_Hs}, 
\begin{equation} \label{eqn:iso}
     \mathcal{L}^\beta: H^q(\Gamma) \rightarrow H^{q-2\beta}(\Gamma) 
\end{equation}
is an isomorphism provided that
$$q\in \lbrack -1-\alpha_{\Gamma}+2\beta,1+\alpha_{\Gamma}\rbrack
\quad \text{for } \beta\in (0,1),
$$
and 
$$q\in \lbrack -1-\alpha_{\Gamma},1+\alpha_{\Gamma}+2\beta \rbrack \quad  \text{for } \beta\in (-1,0).
$$
Moreover, in view of \eqref{def:spaceHs_innerproduct}, we have for $q\in [-1-\alpha_\Gamma,1+\alpha_\Gamma]$
\begin{equation} \label{eq:IP_Hs}
\langle v,w\rangle_{\dot H^q(\Gamma)}=\langle \mathcal{L}^{\frac{q}{2}}v,\mathcal{L}^{\frac{q}{2}}w\rangle_{L^2(\Gamma)}=\langle \mathcal{L}^q v,w\rangle, \qquad v,w\in   H^q(\Gamma)
\end{equation}
and 
\begin{equation} \label{eq:Norm_Hs}
\|v\|_{\dot H^q(\Gamma)}=\|\mathcal{L}^{\frac{q}{2}}v\|_{L^2(\Gamma)}, \qquad v\in  H^q(\Gamma).
\end{equation}

We equipp $\mathcal H^s(\Omega)$ with the $\dot H^s(\Gamma)$ inner product and define the Riesz representer $\dot \phi^i \in \mathcal H^s(\Omega)$ of the functional $\lambda^i:\mathcal H^s(\Omega)\rightarrow\mathbb R$ by the relation
\begin{equation}\label{e:riesz_dot}
\langle \dot \phi^i,v \rangle_{\dot H^s(\Gamma)} = \lambda^i(v), \qquad \forall\, v \in \mathcal H^s(\Omega).
\end{equation}
Since $\dot \phi^i$ is harmonic, it reads $\dot \phi^i = E \dot \psi^i$,
where $E$ is the harmonic extension operator defined for $g\in H^{1/2}(\Gamma)$ as the unique function $w=Eg \in H^1(\Omega)$ satisfying
\begin{equation} \label{eq:E}
w|_{\Gamma} =g \qquad \text{and} \qquad \int_{\Omega}\nabla w \cdot\nabla v = 0, \qquad \forall\, v\in H_0^1(\Omega).
\end{equation}
It is well-known that $E:H^{1/2}(\Gamma) \rightarrow H^1(\Omega)$ is continuous and we denote by $C_E$ the positive constant such that
\begin{equation}\label{e:stab_E}
\| E(v) \|_{H^1(\Omega)} \leq C_E \| v \|_{H^{1/2}(\Gamma)}, \qquad \forall\, v \in H^{1/2}(\Gamma).
\end{equation}

Returning to \eqref{e:riesz_dot} and recalling \eqref{eq:IP_Hs}, we find that the trace $\dot \psi^i \in H^s(\Gamma)$ of $\dot \phi^i$ satisfies
\begin{equation}\label{e:riesz_bdy_dot}
\langle \mathcal L^{s}\dot \psi^i, \eta\rangle = (\lambda^i \circ E)(\eta), \qquad \forall\, \eta \in H^s(\Gamma).
\end{equation}

Note that $\lambda^i \in H^1(\Omega)^*$ implies that 
\begin{equation}\label{e:reg_mu}
\mu^i:=\lambda^i\circ E \in H^{-1/2}(\Gamma)
\end{equation}
with
\begin{align}
     \| \mu^i \|_{H^{-1/2}(\Gamma)} &= \sup_{v \in H^{1/2}(\Gamma) \setminus \{0\}} \frac{\lambda^i(E(v))}{\| v \|_{H^{1/2}(\Gamma)}} 
\leq \| \lambda^i \|_{H^1(\Omega)^*} \sup_{v \in H^{1/2}(\Gamma) \setminus \{0\}}  \frac{\| E(v)\|_{H^1(\Omega)}}{\| v \|_{H^{1/2}(\Gamma)} }  \nonumber \\
& \leq C_E  \| \lambda^i \|_{H^1(\Omega)^*}. \label{eqn:reg_mu}    
\end{align}
Whence, $\dot \psi^i = \mathcal L^{-s}\mu^i \in H^{\min(2s-\frac12,1+\alpha_\Gamma)}(\Gamma)\subset H^s(\Gamma)$, $1/2<s<1$, and the Balakrishnan formula \cite{B1960}, see also Proposition 3.1 in \cite{BGL2024}, offers an integral representation of $\dot \psi^i$
\begin{align} \label{eq:psi_Balak}
\dot \psi^i=\mathcal{L}^{-s}\mu^i &= \frac{\sin(\pi s)}{\pi}\int_{-\infty}^{\infty}e^{(1-s)y}(e^yI+\mathcal{L})^{-1}\mu^i{\rm d}y \nonumber \\
 &= \frac{\sin(\pi s)}{\pi}\int_{-\infty}^{\infty}e^{(1-s)y}w(y){\rm d}y,
\end{align}
where $w(y)\in H^1(\Gamma)$ solves
\begin{equation} \label{eq:w_y}
(e^y+1)\int_{\Gamma}w(y)v+\int_{\Gamma}\nabla_{\Gamma}w(y)\cdot\nabla_{\Gamma}v = \mu^i(v), \qquad \forall\, v\in H^1(\Gamma).
\end{equation}

\section{Computation of the Riesz representers}
\label{sec:algos}

The key ingredient in the recovery algorithm (Algorithm~\ref{alg:OR}) is the approximation for $i=1,...,m$ of the Riesz representer $\dot \phi^i$ via its trace $\dot \psi^i$ satisfying \eqref{eq:psi_Balak}. 
We present two algorithms serving different purposes.
The \emph{theoretical algorithm} is based on the representation \eqref{e:riesz_bdy_dot} and  assumes as in \cite{BBCDDP2023} that an exact computation of the  nonlocal $\dot H^s(\Gamma)$ inner product can be performed. Instead, the \emph{practical algorithm} build from the theoretical algorithm and accounts for the approximation of the $\dot\psi^i$ without computing $\dot H^s(\Gamma)$ inner products. Before describing the two algorithms, we start with the introduction of the finite element tools. For the rest of this section, we simplify the notation and drop the superscript $i$ indexing the measurements.

\subsection{Finite element tools}\label{sec:tools}

Let $\{\mathcal{T}_h\}_{h>0}$ be a sequence of quasi-uniform and shape-regular subdivisions of $\overline{\Omega}$ made of simplices (resp. quadrilaterals/hexahedrons). Associated to each subdivision $\mathcal T_h$, we denote by $\{\varphi_1,\ldots,\varphi_N\}$ and $\{\varphi_{N+1},\ldots,\varphi_{N+N_b}\}$ the Lagrange nodal basis functions \cite{brenner2008mathematical} associated with interior and boundary nodes, respectively.
We define
\begin{equation}\label{e:Vh}
\mathbb{V}_h := \mathbb W_h \oplus \mathbb V_{b,h} \subset H^1(\Omega),
\end{equation}
where
\begin{align*}  
  \mathbb{W}_h &:= {\rm span}\{\varphi_1,\varphi_2,\ldots,\varphi_N\}\subset H_0^1(\Omega) \\
  \mathbb{V}_{b,h} & := {\rm span}\{\varphi_{N+1},\varphi_{N+2},\ldots,\varphi_{N+N_b}\}\subset H^1(\Omega)
\end{align*}
and
$$
\mathbb{T}_h  := \mathbb{V}_{b,h}|_{\Gamma}=\mathbb{V}_{h}|_{\Gamma}={\rm span}\{\varphi_{N+1}|_{\Gamma},\varphi_{N+2}|_{\Gamma},\ldots,\varphi_{N+N_b}|_{\Gamma}\}\subset H^s(\Gamma),
$$ 
since $\varphi_i|_\Gamma=0$ for any $i\in\{1,2,\ldots,N\}$. 

We can now define the discrete harmonic extension operator $E_h:\mathbb T_h \rightarrow \mathbb V_h$. For $g_h \in \mathbb T_h$, $w_h:=E_h g_h \in \mathbb V_h$ is given by 
\begin{equation} \label{eq:E_h}
w_h |_\Gamma=g_h \qquad \text{and} \qquad \int_{\Omega}\nabla w_h \cdot\nabla v_h = 0 \quad \forall\, v_h\in \mathbb{W}_h;
\end{equation}
compare with \eqref{eq:E}. 
The stability of $E_h$ follows from the stability \eqref{e:stab_E} of $E$, refer to \cite{BBCDDP2023} for more details.
In particular, there is a constant $D_E>0$ independent of $h$ such that 
\begin{equation} \label{eqn:stab_E_h}
    \|E_hg_h\|_{H^1(\Omega)}\le D_E\|g_h\|_{H^{1/2}(\Gamma)}, \qquad \forall\, g_h\in\mathbb{T}_h.
\end{equation}

We shall also need the Scott-Zhang projection $\Pi_h : L^2(\Gamma) \rightarrow \mathbb V_h$ introduced in \cite{scott1990finite} for Euclidean domains but whose definition and properties readily extend to the boundary of a polygonal domain like $\Gamma$. Its stability in $L^2(\Gamma)$ and $H^1(\Gamma)$ are standard and for $t\in [0,1]$ we have
\begin{equation} \label{eq:stab_pi_h}
    \|\Pi_h v\|_{H^{t}(\Gamma)}\lesssim \|v\|_{H^t(\Gamma)}, \qquad \forall\, v\in H^t(\Gamma).
\end{equation}
In \cite{faustmann2021stability} the stability property \eqref{eq:stab_pi_h} is extended to $t\in (1,3/2)$ but with the Hilbert norms replaced by Besov norms. Nevertheless, we have for $t\in (1,3/2)$
\begin{equation} \label{eq:stab_pi_h_more}
    \|\Pi_h v\|_{H^{t-\epsilon}(\Gamma)}\lesssim C(\epsilon) \|v\|_{H^t(\Gamma)}, \qquad \forall\, v\in H^t(\Gamma),
\end{equation}
where $\epsilon>0$ and $C(\epsilon) \to \infty$ as $\epsilon \to 0^+$.
Moreover, for any $t_1 \in [0,1]$, $t_2\in[0,2]$ such that $t_1+t_2 \in [0,2]$, we have 
\begin{equation} \label{eq:est_pi_h}
\|(I-\Pi_h) v\|_{H^{t_1}(\Gamma)}\lesssim h^{t_2}\|v\|_{H^{t_1+t_2}(\Gamma)}, \qquad \forall\, v\in H^{t_1+t_2}(\Gamma).
\end{equation}

We shall also need the discrete counter-part of the operator $\mathcal L:=(I-\Delta_\Gamma)$, namely $\mathcal L_h: \mathbb T_h \rightarrow \mathbb T_h$, where for $w_h \in \mathbb T_h$, $g_h=\mathcal L_h w_h$ satisfies
$$
\int_\Gamma g_h v_h =  \int_\Gamma w_h v_h + \int_\Gamma \nabla_\Gamma w_h \cdot \nabla_\Gamma v_h, \qquad \forall\, v_h \in \mathbb T_h.
$$

\subsection{Theoretical algorithm} \label{sec:theoretical}

The finite element approximation of the Riesz representer $\dot \phi\in H^1(\Omega)$ is defined by $\dot \phi_h := E_h \dot \psi_h\in \mathbb{V}_h$, where $\dot \psi_h\in \mathbb{T}_h$ approximates $\dot \psi$ and is given by the relations
\begin{equation} \label{eq:Riesz_h}
    \langle \dot \psi_h,\eta_h\rangle_{\dot H^s(\Gamma)} = \mu_h(\eta_h):=(\lambda \circ E_h)(\eta_h), \qquad \forall\, \eta_h\in \mathbb{T}_h.
\end{equation}
We have already pointed out that the computation of $\dot\phi_h=E_h\dot\psi_h$ using this approach is quite demanding for two reasons: for all $1\le i,j \le N_b$ one must compute
\begin{enumerate}
    \item exactly the nonlocal $H^s(\Gamma)$ inner product $\langle \varphi_{N+j},\varphi_{N+i}\rangle_{\dot H^s(\Gamma)}$ leading to a full system; 
    \item the discrete harmonic extension $E_h (\varphi_{N+i}|_\Gamma)$.
    \end{enumerate}

The first issue is addressed in the next section. 
For the second issue, we follow \cite{BBCDDP2023} and consider the following saddle-point problem enforcing the harmonic constraint using Lagrange multipliers:
find $(\dot \phi_h,\xi_h)\in \mathbb{V}_h\times\mathbb{W}_h$ such that
\begin{equation} \label{eq:saddle_h}
\begin{array}{lcll}
\langle \dot\phi_h,v_h\rangle_{\dot H^s(\Gamma)}+\int_\Omega \nabla v_h \cdot \nabla \xi_h & = & \lambda(v_h) &\qquad  \forall\, v_h\in \mathbb{V}_h, \medskip\\
\int_\Omega \nabla \dot \phi_h \cdot \nabla z_h & = & 0 & \qquad  \forall\, z_h\in \mathbb{W}_h.
\end{array}
\end{equation}
In \cite{BBCDDP2023}, it is shown that the above saddle-point problem has a unique solution $(\dot \phi_h,\xi_h)$ and  $\dot \phi_h = E_h \dot \psi_h$, where $\dot \psi_h \in \mathbb T_h$ solves \eqref{eq:Riesz_h}.   
We now further explore the structure of \eqref{eq:saddle_h} to reduce its complexity
and write $\dot \phi_h = \dot \phi_h^0 + \dot \phi_h^b$, where $\dot \phi_h^0 \in \mathbb W_h$ and $\dot \phi_h^b \in \mathbb V_{b,h}$. 
Restricting the first equation of \eqref{eq:saddle_h} to $v_h \in \mathbb W_h \subset \mathbb V_h$  yields 
\begin{equation}\label{e:pi_h_saddle}
\int_\Omega \nabla \xi_h \cdot \nabla v_h  = \lambda(v_h), \qquad \forall\, v_h \in \mathbb W_h,
\end{equation}
which uniquely determine $\xi_h \in \mathbb W_h$.
Then, restricting the same equation to  $v_h \in \mathbb V_{b,h}\subset \mathbb V_h$ implies
\begin{equation}\label{e:phi_bh_saddle}
\langle \dot\phi_h^b ,v_h\rangle_{\dot H^s(\Gamma)} =  \lambda(v_h) - \int_\Omega \nabla \xi_h \cdot \nabla v_h, \qquad \forall\, v_h \in \mathbb V_{b,h},
\end{equation}
which this time uniquely determine $\dot \phi_h^b \in \mathbb V_{b,h}$ and thus $\dot \phi_{h}^b|_{\Gamma} \in \mathbb T_h$.
It remains to use the second equation in \eqref{eq:saddle_h} to compute $\dot \phi_h^0 \in \mathbb W_h$ solving
\begin{equation}\label{e:phi_0h_saddle}
\int_{\Omega} \nabla \dot \phi_h^0 \cdot \nabla v_h = - \int_\Omega \nabla \phi_h^b \cdot \nabla v_h, \qquad \forall\, v_h \in \mathbb W_h.
\end{equation}

We summarize this sequential algorithm in Algorithm~\ref{alg:theoretical}.
\begin{algorithm}
\caption{Theoretical Algorithm for the approximation of $\dot \phi$}\label{alg:theoretical}
\begin{algorithmic}
\State $\diamond$ Find $\xi_h \in \mathbb W_h$ satisfying \eqref{e:pi_h_saddle};
\State $\diamond$ Find $\dot \phi_h^b \in \mathbb V_{b,h}$ satisfying \eqref{e:phi_bh_saddle};
\State $\diamond$ Find $\dot \phi_h^0 \in \mathbb W_{h}$ satisfying \eqref{e:phi_0h_saddle};
\State $\diamond$ Compute $\dot \phi_h = \dot \phi_h^0 + \dot \phi_h^b$.
\end{algorithmic}
\end{algorithm}

We end this section by noting that the sequential resolution of \eqref{e:pi_h_saddle}, \eqref{e:phi_bh_saddle}, and \eqref{e:phi_0h_saddle} has the advantage of not requiring to assemble and store the (potentially large) matrix of size $(2N+N_b)\times(2N+N_b)$ associated to the saddle-point problem \eqref{eq:saddle_h}. More importantly, in contrast to the formulation \eqref{eq:Riesz_h}, it does not require the expensive computation of the discrete harmonic extensions $E_h(\varphi_{N+i}|_\Gamma)$ for $i=1,2,\ldots,N_b$. 
However, it still assumes the exact computation of $H^s(\Gamma)$ inner products $\langle \varphi_{N+j},\varphi_{N+i}\rangle_{\dot H^s(\Gamma)}$ to assemble the matrix associated to Problem \eqref{e:phi_bh_saddle}. We address this aspect in the next section.

\subsection{Practical algorithm} \label{sec:practical}

We restart from the Balakrishnan formula \eqref{eq:psi_Balak} for $\dot \psi$ and proceed as in \cite{BLP2019,BL2022} for its approximation, thereby circumventing the computation of the $\dot H^s(\Gamma)$ inner product in \eqref{eq:Riesz_h}. A sinc quadrature approximate the indefinite integral and a finite element method is advocated for the approximation of $w(y_l)$ at each quadrature point $y_l$. 
 
 Given $k>0$, we set $y_l:=lk$ for $l=-\texttt{M},\ldots,\texttt{N}$, where the specific values for $\texttt{M},\texttt{N}\in\mathbb{N}$ will be chosen later (see Lemma~\ref{prop:Sinc}) to guarantee a quadrature exhibiting an exponential rate of approximation. The sinc approximation of $\dot \psi$ reads
\begin{align} \label{eq:psi_sinc}
\dot \psi_k := \mathcal{Q}_k^{-s}(\mathcal L)\mu &:= \frac{k\sin(\pi s)}{\pi}\sum_{l=-\texttt{M}}^\texttt{N}e^{(1-s)y_l}(e^{y_l}I+\mathcal{L})^{-1}\mu \nonumber \\
&= \frac{k\sin(\pi s)}{\pi}\sum_{l=-\texttt{M}}^\texttt{N}e^{(1-s)y_l}w(y_l),
\end{align}
where $w(y_l)\in H^1(\Gamma)$ solves \eqref{eq:w_y}.
We  further approximate $\dot\psi_k$ using a finite element method based on $\mathbb T_h$, namely we define
\begin{equation} \label{eq:psi_sinc_FEM}
\dot \psi_{k,h}:= \frac{k\sin(\pi s)}{\pi}\sum_{l=-\texttt{M}}^\texttt{N}e^{(1-s)y_l}w_h(y_l),
\end{equation}
where for $y\in\mathbb{R}$, $w_h(y)\in\mathbb{T}_h$ is the solution to
\begin{equation} \label{eq:w_y_FEM}
(e^y+1)\int_{\Gamma}w_h(y)v_h+\int_{\Gamma}\nabla_{\Gamma}w_h(y)\cdot\nabla_{\Gamma}v_h = \mu_h(v_h), \qquad \forall\, v_h\in\mathbb{T}_h,
\end{equation}
where we recall that $\mu_h(v_h)=\lambda(E_h v_h)$ for $v_h \in \mathbb T_h$. 
Using the operator $\mathcal L_h$ defined in Section~\ref{sec:tools} and recalling that $L^2(\Gamma)$ is identified with its dual so that $\mu_h(v_h) = \int_\Gamma \mu_h v_h$ for all $v_h \in \mathbb T_h$, we have
$$
\dot \psi_{k,h}= \mathcal{Q}_k^{-s}(\mathcal L_h)\mu_h := \frac{k\sin(\pi s)}{\pi}\sum_{l=-\texttt{M}}^\texttt{N}e^{(1-s)y_l} (e^{y_l}I+\mathcal L_h)^{-1}\mu_h.
$$

Once the approximation $\dot \psi_{k,h}$ is computed, the approximation of the Riesz representer $\dot \phi$ is then given by $\dot \phi_{k,h}:=E_h \dot \psi_{k,h}$. In Section~\ref{sec:error}, we analyze the behavior of $\| \dot\phi - \dot\phi_{k,h} \|_{H^1(\Omega)}$. At this point, we emphasize again that \eqref{eq:psi_sinc_FEM} does not require the computation of $H^s(\Gamma)$ inner products; compare with \eqref{eq:Riesz_h}. 
Also, notice that $\dot \psi_{k,h}$ is obtained by aggregating $\texttt{N}+\texttt{M}+1$ standard approximation of reaction-diffusion problems.

As for the theoretical algorithm (Algorithm~\ref{alg:theoretical}), we circumvent the computation of the $N_b$ discrete harmonic extensions using the saddle-point framework.
We write $\dot \phi_{k,h} = \dot \phi_{k,h}^0 + \dot \phi_{k,h}^b$ with $\dot \phi_{k,h}^0\in \mathbb W_h$ and $\dot \phi_{k,h}^b \in \mathbb V_{b,h}$. 
Step 1 and 3 in Algorithm~\ref{alg:theoretical} remain unchanged while step 2 is replaced by 
\begin{equation}\label{eq:psi_sinc_saddle}
\dot \phi_{k,h}^b|_\Gamma = \frac{k\sin(\pi s)}{\pi}\sum_{l=-\texttt{M}}^\texttt{N}e^{(1-s)y_l}w_h(y_l),
\end{equation}
where for $y\in\mathbb{R}$, $w_h(y) \in \mathbb T_h$ is given by
\begin{equation} \label{eq:w_y_FEM_sinc}
(e^y+1)\int_{\Gamma}w_h(y)v_h+\int_{\Gamma}\nabla_{\Gamma}w_h(y)\cdot\nabla_{\Gamma}v_h = \lambda(v_h) - \int_\Omega \nabla \xi_h \cdot \nabla v_h, \quad \forall\, v_h\in\mathbb{V}_{b,h},
\end{equation}
instead of \eqref{eq:w_y_FEM}. The full practical algorithm is described in Algorithm~\ref{alg:practical}.
\begin{algorithm}
\caption{Practical Algorithm for the approximations of $\dot \phi$}\label{alg:practical}
\begin{algorithmic}
\State $\diamond$ Find $\xi_h \in \mathbb W_h$ satisfying \eqref{e:pi_h_saddle};
\State $\diamond$ Find $\dot \phi_{k,h}^b \in \mathbb V_{b,h}$ satisfying \eqref{eq:psi_sinc_saddle};
\State $\diamond$ Find $\dot \phi_h^0 \in \mathbb W_{h}$ satisfying \eqref{e:phi_0h_saddle} with $ \phi_h^b$ replaced by $\dot \phi_{k,h}^b$;
\State $\diamond$ Compute $\dot \phi_{k,h} = \dot \phi_h^0 + \dot \phi_{k,h}^b$.
\end{algorithmic}
\end{algorithm}

The next lemma guarantees that the output $\dot \phi_{k,h}$ of Algorithm~\ref{alg:practical} is such that $\dot \phi_{k,h} = E_h \dot \psi_{k,h}$, where 
$\dot \psi_{k,h} \in \mathbb T_h$ satisfies \eqref{eq:psi_sinc_FEM}.
\begin{lemma}
Algorithm~\ref{alg:practical} produces the unique discrete harmonic function $\dot \phi_{k,h} \in \mathbb V_h$ whose trace  $\dot \phi_{k,h}|_\Gamma$ satisfies \eqref{eq:psi_sinc_FEM}.
\end{lemma}
\begin{proof}
    Comparing \eqref{eq:psi_sinc_FEM}-\eqref{eq:w_y_FEM} with \eqref{eq:psi_sinc_saddle}-\eqref{eq:w_y_FEM_sinc}, it suffices to show that
\begin{equation}\label{e:same_rhs}
\lambda(v_h) - \int_\Omega \nabla \xi_h \cdot \nabla v_h = \lambda (E_h (v_h|_\Gamma))=\mu_h(v_h|_\Gamma)
\end{equation}
for all $v_h \in \mathbb V_{b,h}$. 
Since $v_h - E_h (v_h|_\Gamma) \in \mathbb W_h$,
relation \eqref{e:pi_h_saddle} together with the definition of the discrete harmonic extension operator \eqref{eq:E_h} yield
$$
\lambda(v_h - E_h(v_h|_\Gamma)) = \int_\Omega \nabla \xi_h \cdot \nabla (v_h - E_h(v_h|_\Gamma)) = \int_\Omega \nabla \xi_h \cdot \nabla v_h,
$$
which is \eqref{e:same_rhs} in disguised. From \eqref{e:same_rhs} we deduce that 
$\dot \phi_{b,h}|_\Gamma = \dot \psi_h$, where $\dot\psi_h \in \mathbb T_h$ solves \eqref{eq:Riesz_h}. The proof is complete.
\end{proof}

The representation \eqref{eq:psi_sinc_FEM} is instrumental in the analysis provided below while Algorithm~\ref{alg:practical} is advocated in Section~\ref{sec:num_res} to illustrate the performance of the method.

\section{Error analysis for the Riesz representer and optimal recovery}
\label{sec:error}

As in the previous section, we drop the notation $i$ indexing the measurements and let $\lambda \in H^{1}(\Omega)^*$ be a generic measurement with associated $\mu:=\lambda \circ E$ which belongs to $H^{-1/2}(\Gamma)$ in view of \eqref{e:reg_mu}.

The efficiency of the recovery algorithm (Algorithm~\ref{alg:OR}) depends on the discrepancy between $\dot \phi=E\dot \psi$ with $\dot \psi$ the solution to the boundary problem \eqref{e:riesz_bdy_dot} and $\dot \phi_{k,h}=E_h\dot \psi_{k,h}$ with $\dot \psi_{k,h}$ is the sinc-FE approximation defined in \eqref{eq:psi_sinc_FEM}.

In this section, we derive an error estimate for $\| \dot \phi - \dot \phi_{k,h} \|_{H^1(\Omega)}$ by decomposing the error into three parts
\begin{equation} \label{eqn:splitting_error}
\|\dot\phi-\dot\phi_{k,h}\|_{H^1(\Omega)} = \|E\dot\psi-E_h\dot\psi_{k,h}\|_{H^1(\Omega)} \le {\rm I}+{\rm II}+{\rm III},
\end{equation}
where
\begin{equation*}
    \begin{split}
    &{\rm I} :=\|E\dot\psi-E\dot\psi_k\|_{H^1(\Omega)}, \qquad {\rm II} :=\|E\dot\psi_k-E_h(\Pi_h\dot\psi_k)\|_{H^1(\Omega)}, \\
    &\textrm{and} \qquad {\rm III} :=\|E_h(\Pi_h\dot\psi_k)-E_h\dot\psi_{k,h}\|_{H^1(\Omega)}
    \end{split}
\end{equation*}
correspond to the sinc quadrature error, the finite element error for the discrete harmonic extension, and the finite element error for $\dot \psi_k$, respectively. 
We recall that $\Pi_h$ denotes the Scott-Zhang operator defined in Section~\ref{sec:tools} and that we restrict our attention to the case $s\in(1/2,1)$.  

We bound each term separately in the following three subsections to derive the following result.

\begin{theorem}\label{t:main}
      Let $1/2<s<1$, $\alpha_\Gamma$ be as in Assumption~\ref{a:pickup} and set
      \begin{equation} \label{eqn:def_gamma_tilde}
      \widetilde\gamma:=\min(2s-1,\alpha_\Gamma+1/2,2\alpha_\Gamma).  
    \end{equation}
    Let $\alpha_\Omega$ be as in Assumption~\ref{ass:E}. For $i=1,...,m$, let $\lambda^i \in H^1(\Omega)^*$, $\dot \phi^i$ be given by \eqref{e:riesz_dot} and $\dot \phi^i_{k,h}$ be the output of Algorithm~\ref{alg:practical}. For every $0<\epsilon<\widetilde\gamma$, there are constants $C(\epsilon)$ and $C$ such that
\begin{equation}
\| \dot \phi^i - \dot \phi^i_{k,h} \|_{{H^{1}(\Omega)}} \leq C(\epsilon) h^{\min\{\widetilde\gamma-\epsilon,\alpha_{\Omega}\}}+ C e^{-\frac{\pi^2}{k}}\|\lambda^i\|_{H^{1}(\Omega)^*},
\end{equation}
where $C(\epsilon) \to \infty$ as $\epsilon \to 0^+$.
\end{theorem}
\begin{proof}
The proof follows upon combining Corollary~\ref{c:sinc}, Lemma~\ref{lem:boundII} and Lemma~\ref{lem:boundIII} below.
\end{proof}
\begin{corollary}\label{c:main}
   Under the assumptions of Theorem~\ref{t:main}, assume further that the model class Assumption~\ref{ass:model} holds. Fix $0<\epsilon<\widetilde\gamma$ and let 
    $$
    \varepsilon = \left(C(\epsilon) h^{\min\{\widetilde\gamma-\epsilon,\alpha_{\Omega}\}}+ C e^{-\frac{\pi^2}{k}}\right)\max_{i=1,...,m}\|\lambda^i\|_{H^{1}(\Omega)^*},
    $$
    where $C(\epsilon)$ and $C$ are the constants in Theorem~\ref{t:main}.
    Then the output $\hat u = u_{k,h}$ of Algorithm~\ref{alg:OR} with input parameter $\varepsilon$ and where the Riesz representers are computed using Algorithm~\ref{alg:practical} satisfies the near optimal recovery estimate 
    $$
    \| u- u_{k,h} \|_{H^1(\Omega)} \leq R(\dot K^s_{\boldsymbol \omega}) + \dot \Lambda \left(C(\epsilon) h^{\min\{\widetilde\gamma-\epsilon,\alpha_{\Omega}\}}+ C e^{-\frac{\pi^2}{k}})\right)\max_{i=1,...,m}\|\lambda^i\|_{H^{1}(\Omega)^*},
    $$
    where $\dot \Lambda$ is the constant appearing in \eqref{e:OR_dot}.
\end{corollary}
\begin{proof}
Recall that we assume that $u_f=f=0$ so it suffices to combine the optimal recovery estimate \eqref{e:OR_dot} with the specific choice of $\varepsilon$ and Theorem~\ref{t:main}.
\end{proof}
\subsection{The sinc quadrature error} \label{est:boundI} %

We derive in this section an estimate for the sinc quadrature error, i.e., for ${\rm I}$ in \eqref{eqn:splitting_error}. We take advantage of the stability \eqref{e:stab_E} of the harmonic extension operator $E$ and the following lemma.


\begin{lemma} \label{prop:Sinc}
Let  $\lambda \in H^1(\Omega)^*$, $\dot \psi=\mathcal{L}^{-s}(\lambda \circ E)$ be as in \eqref{eq:psi_Balak} and $\dot \psi_k=\mathcal{Q}_k^{-s}(\mathcal{L})(\lambda \circ E)$ be as in \eqref{eq:psi_sinc} with
\begin{equation} \label{def:M_N_sinc}
    \emph{\texttt{M}}:=\left\lceil\frac{\pi^2}{(1-s)k^2}\right\rceil \qquad \text{and} \qquad \emph{\texttt{N}}:=\left\lceil\frac{\pi^2}{\left(s-\frac{1}{2}\right)k^2}\right\rceil.
\end{equation}
Then we have
$$\|\dot\psi-\dot\psi_k\|_{H^{1/2}(\Gamma)}\lesssim e^{-\frac{\pi^2}{k}}\| \lambda\|_{H^{1}(\Omega)^*},$$
where the hidden constant depends on $\max\{(1-s)^{-1},(2s-1)^{-1}\}$.
\end{lemma}
\begin{proof}
We only give a sketch of the proof as it follows from \cite[Theorem 3.2]{BLP2019} (case $\lambda \circ E \in L^2(\Gamma)$) and \cite[Proposition 5.1]{BGL2024} (estimate for $\| \dot \psi - \dot \psi_k \|_{L^2(\Gamma)}$).

We set $\mu:= \lambda \circ E$, which belong to $H^{-1/2}(\Gamma) \cong \dot H^{-1/2}(\Gamma)$ according to \eqref{e:reg_mu} and \eqref{e:equiv_Hs}. Whence, we have $\dot \psi=\mathcal{L}^{-s}\mu\in \dot H^{2s-1/2}(\Gamma)$. 
Furthermore, we write
$$\|\dot \psi-\dot\psi_{k}\|_{\dot H^{1/2}(\Gamma)}=\|\mathcal{L}^{\frac{1}{4}}(\dot\psi-\dot \psi_{k})\|_{L^2(\Gamma)}=\|\mathcal{L}^{\frac{1}{4}}(\mathcal{L}^{-s}-\mathcal{Q}_k^{-s}(\mathcal{L}))\mathcal{L}^{\frac{1}{4}}v\|_{L^2(\Gamma)},$$
where $v:=\mathcal{L}^{-1/4}\mu\in L^2(\Gamma)$.
In view of the representation \eqref{eq:psi_Balak} of $\dot\psi$ and \eqref{eq:psi_sinc} of $\dot\psi_k$, we expand $v$ in the basis of $L^2(\Gamma)$-orthonormal eigenfunctions $\{e_j\}_{j=1}^{\infty}$ according to \eqref{eqn:eig_expansion} and deduce that for all $y\in\mathbb{R}$
\begin{align*}
    e^{(1-s)y}\mathcal{L}^{\frac{1}{4}}(e^yI+\mathcal{L})^{-1}\mathcal{L}^{\frac{1}{4}}v &=e^{(1-s)y}\sum_{j=1}^{\infty}v_j\mathcal{L}^{\frac{1}{4}}(e^yI+\mathcal{L})^{-1}\mathcal{L}^{\frac{1}{4}}e_j \\
    &=e^{(1-s)y}\sum_{j=1}^{\infty}v_j\eig_j^{\frac{1}{4}}(e^y+\eig_j)^{-1}\eig_j^{\frac{1}{4}}e_j,
\end{align*}
where $v_j:= \langle v, e_j \rangle_{L^2(\Gamma)}$.
This implies that
$$\|\dot\psi-\dot\psi_{k}\|_{\dot H^{1/2}(\Gamma)}\le \left(\max_{\eig\ge \eig_1}|e(\eig)|\right)\|v\|_{L^2(\Gamma)}=\left(\max_{\eig\ge \eig_1}|e(\eig)|\right)\|\mu\|_{H^{-1/2}(\Gamma)}$$
with
$$e(\eig):=\int_{-\infty}^{\infty}g_{\eig}(y){\rm d}y-k\sum_{l=-\texttt{M}}^{\texttt{N}}g_{\eig}(y_l), \quad g_{\eig}(y):=e^{(1-s)y}(e^y+\eig)^{-1}\eig^{\frac{1}{2}}.$$
The choice \eqref{def:M_N_sinc} for $\texttt{M}$ and $\texttt{N}$ together with Lemma 3.1 in \cite{BLP2019}, see also \cite{LB1992}, guarantee that
$$e(\eig)\lesssim e^{-\frac{\pi^2}{k}}+e^{-\left(s-\frac{1}{2}\right)\texttt{N}k}+e^{-(1-s)\texttt{M}k}$$
so that
$$
\|\dot\psi-\dot\psi_{k}\|_{\dot H^{1/2}(\Gamma)} \lesssim \left( e^{-\frac{\pi^2}{k}}+e^{-\left(s-\frac{1}{2}\right)\texttt{N}k}+e^{-(1-s)\texttt{M}k}\right) \| \mu \|_{H^{-1/2}(\Gamma)}.
$$
The desired result follows upon using again the equivalence \eqref{e:equiv_Hs} and estimate \eqref{eqn:reg_mu} for $\| \mu\|_{H^{-1/2}(\Gamma)}$.
\end{proof}
\begin{corollary}\label{c:sinc}
    Under the assumptions of Lemma~\ref{prop:Sinc}, there holds
    \begin{equation} \label{eqn:bound_I}
     \|E\psi-E\psi_k\|_{H^1(\Omega)}\lesssim e^{-\frac{\pi^2}{k}}\|\lambda\|_{H^{1}(\Omega)^*}.  
\end{equation}
\end{corollary}
\begin{proof}
This result directly follows from the stability \eqref{e:stab_E} of the extension operator and Lemma~\ref{prop:Sinc}.
\end{proof}

\subsection{The finite element error for the harmonic extension}
\label{sec:E_Eh}

In this section, we derive an upper bound for the term ${\rm II}$ in \eqref{eqn:splitting_error}. We start by estimating the discrepancy between $Eg_h$ and $E_h g_h$ for $g_h \in \mathbb T_h$, which depends on the regularity of $Eg_h$ for $g_h \in \mathbb T_h$. We make the following assumption.
\begin{assumption}\label{ass:E}
    There is $\alpha_\Omega>0$ and a constant $C$ such that for all $q \in [0,\alpha_\Omega]$ there holds
    $$
    Ev \in H^{1+q}(\Omega) \qquad \textrm{with} \qquad 
    \| Ev \|_{H^{1+q}(\Omega)} \leq C \| v \|_{H^{1/2+q}(\Gamma)},
    $$
    for all $v \in H^{1/2+q}(\Gamma)$.
\end{assumption}
As suggested by the notation, the value of $\alpha_\Omega$ depends on $\Omega$. For instance, for polygonal domains $\Omega$ with inner angles at most $\theta$, we have 
that $\alpha_{\Omega}=\pi/\theta-\epsilon$ for all $\epsilon>0$; see \cite{G1985}.

The approximation of the harmonic extension operator $E$ by its discrete counter-part $E_h$ is now assessed under Assumption~\ref{ass:E}.

\begin{lemma} \label{prop:extensionE}
    Let $\alpha_\Omega$ be as in Assumption~\ref{ass:E}. For any $g_h\in\mathbb{T}_h$ and any $\beta\in(1/2,3/2)$ we have
    $$\|Eg_h-E_hg_h\|_{H^1(\Omega)}\lesssim h^{\min\{\beta-1/2,\alpha_{\Omega}\}}\|g_h\|_{H^{\beta}(\Gamma)}.$$
\end{lemma}

\begin{proof}
    Lemma 4.3 in \cite{BBCDDP2023} guarantees that
    $$\|Eg_h-E_hg_h\|_{H^1(\Omega)} \lesssim \inf_{v_h\in\mathbb{V}_h}\|Eg_h-v_h\|_{H^1(\Omega)}.$$
    
    Notice that the functions in $\mathbb{T}_h$ belong to $W^{1,\infty}(\Gamma)$ and thus are in $W^{1,p}(\Gamma)$ for any $p\ge 1$. Consequently, we have $g_h\in H^{\beta}(\Gamma)$ for any $\beta<3/2$ by the Sobolev embedding theorem.
    Whence, Assumption~\ref{ass:E} implies that $Eg_h\in H^{1+\min\{\beta-1/2,\alpha_{\Omega}\}}(\Omega)$ and the desired result follows from a standard finite element approximation estimate.
\end{proof}

We now use the interpolation result \eqref{eq:est_pi_h}, the stability estimates \eqref{eq:stab_pi_h}-\eqref{eq:stab_pi_h_more}, and Lemma~\ref{prop:extensionE} to estimate ${\rm II}$.

\begin{lemma} \label{lem:boundII}
    Let $1/2<s<1$, $\alpha_\Gamma$ be as in Assumption~\ref{a:pickup} and $\widetilde\gamma$ be as in \eqref{eqn:def_gamma_tilde}. Let $\alpha_\Omega$ be as in Assumption~\ref{ass:E}. Moreover for $\lambda \in H^1(\Omega)^*$ we set $\mu=\lambda \circ E$ and let $\dot \psi_k=\mathcal{Q}_k^{-s}(\mathcal{L})\mu$  be given by \eqref{eq:psi_sinc} with $\texttt{M}$, $\texttt{N}$ as in \eqref{def:M_N_sinc}. For any $0<\epsilon<\widetilde{\gamma}$, there is a constant $C(\epsilon)$ such that $C(\epsilon)\to \infty$ as $\epsilon \to 0^+$ and    
\begin{equation} \label{eqn:bound_II_case}
   \|E\dot\psi_k-E_h(\Pi_h\dot\psi_k)\|_{H^1(\Omega)} \leq C(\epsilon) h^{\min\{\widetilde{\gamma}-\epsilon,\alpha_{\Omega}\}}\|\lambda\|_{H^{1}(\Omega)^*}.
\end{equation}
\end{lemma}

\begin{proof}
We first write $E\dot\psi_k-E_h(\Pi_h\dot\psi_k) =  E\dot\psi_k-E(\Pi_h\dot\psi_k) + E(\Pi_h\dot\psi_k)-E_h(\Pi_h\dot\psi_k)$ so that in view of  \eqref{e:stab_E} we get
\begin{align}
    {\rm II} &\le \|E\dot\psi_k-E(\Pi_h\dot\psi_k)\|_{H^1(\Omega)} + \|E(\Pi_h\dot\psi_k)-E_h(\Pi_h\dot\psi_k)\|_{H^1(\Omega)} \nonumber \\
    &\le C_E\|\dot\psi_k-\Pi_h\dot\psi_k\|_{H^{1/2}(\Gamma)} + \|(E-E_h)(\Pi_h\dot\psi_k)\|_{H^1(\Omega)}. \label{eqn:II_split}
\end{align}

The approximation properties \eqref{eq:est_pi_h} of the Scott-Zhang projection with $t_1=1/2$ and $t_2=\widetilde{\gamma} - \epsilon$, yield for the first term in \eqref{eqn:II_split}
    \begin{equation} \label{eqn:II_case1a}
    \|\dot\psi_k-\Pi_h\dot\psi_k\|_{H^{1/2}(\Gamma)}\lesssim h^{\widetilde{\gamma}-\epsilon}\|\dot\psi_k\|_{H^{1/2+\widetilde{\gamma} - \epsilon}(\Gamma)}.
    \end{equation}
    For the second term, we invoke Lemma~\ref{prop:extensionE} with $\beta=1/2 + \widetilde{\gamma}-\epsilon<3/2$ together with the stability estimate \eqref{eq:stab_pi_h_more} to write
    \begin{align}
      \|(E-E_h)(\Pi_h\dot\psi_k)\|_{H^1(\Omega)} &\lesssim h^{\min\{\widetilde{\gamma}-\epsilon,\alpha_{\Omega}\}}\|\Pi_h\dot\psi_k\|_{H^{1/2+\widetilde{\gamma}-\epsilon}(\Gamma)} \nonumber \\
      &\lesssim C(\epsilon) h^{\min\{\widetilde{\gamma}-\epsilon,\alpha_{\Omega}\}}\|\dot\psi_k\|_{H^{1/2 + \widetilde{\gamma}-\epsilon/2}(\Gamma)}, \label{eqn:II_case1b}
    \end{align}
    where $C(\epsilon) \to \infty$ as $\epsilon \to 0^+$.
    Inserting \eqref{eqn:II_case1a}, \eqref{eqn:II_case1b}  in \eqref{eqn:II_split} we get
    $$
    {\rm II} \lesssim  C(\epsilon) h^{\min\{\widetilde{\gamma}-\epsilon,\alpha_{\Omega}\}}\|\dot\psi_k\|_{H^{1/2 + \widetilde{\gamma} - \epsilon/2}(\Gamma)}.
    $$
    To end the proof, it remains to note that the stability of $\mathcal{Q}_k^{-s}(\mathcal{L})$ (Lemma~\ref{prop_stability_Qks} with $r=1/2$ and $\beta=1/2+\widetilde\gamma-\epsilon$) implies that $\dot\psi_k=\mathcal{Q}_k^{-s}(\mathcal{L})\mu \in H^{1/2 + \widetilde{\gamma} -\epsilon/2}(\Gamma)$ and 
\begin{equation}\label{e:stb_psi} 
\| \dot\psi_k\|_{H^{1/2 + \widetilde{\gamma} -\epsilon/2}(\Gamma)} \lesssim C(\epsilon) \| \mu \|_{H^{-1/2}(\Gamma)} \lesssim  C(\epsilon) \| \lambda \|_{H^1(\Omega)^*}.
\end{equation}

\end{proof}

\subsection{The finite element for $\dot \psi_k$} \label{est:boundII}

The following lemma provides an estimate for the last term ${\rm III}$ in \eqref{eqn:splitting_error}.

\begin{lemma} \label{lem:boundIII}
    Let $1/2<s<1$, $\alpha_\Gamma$ be as in Assumption~\ref{a:pickup} and $\widetilde\gamma$ be as in \eqref{eqn:def_gamma_tilde}.
    Let $\alpha_\Omega$ be as in Assumption~\ref{ass:E}. For $\lambda \in H^1(\Omega)^*$, we set $\mu:= \lambda \circ E$ and $\mu_h := \lambda \circ E_h$. Let $\dot \psi_k=\mathcal{Q}_k^{-s}(\mathcal{L})\mu$ be given by \eqref{eq:psi_sinc} with $\texttt{M}$ and $\texttt{N}$ as in \eqref{def:M_N_sinc} and $\dot\psi_{k,h} = \mathcal{Q}_k^{-s}(\mathcal{L}_h)\mu_h$ be as in \eqref{eq:psi_sinc_FEM}.  
    For any $0<\epsilon<\widetilde\gamma$, there is a constant $C(\epsilon)$ such that $C(\epsilon)\to \infty$ as $\epsilon \to 0^+$ and    
\begin{equation} \label{eqn:bound_III_case}
   \| E_h(\Pi_h\dot\psi_k)-E_h\dot\psi_{k,h}\|_{H^1(\Omega)} \leq C(\epsilon) h^{\min\{\widetilde\gamma-\epsilon,\alpha_{\Omega}\}}\|\lambda\|_{H^{1}(\Omega)^*}.
\end{equation}
\end{lemma}
\begin{proof}
From the stability \eqref{eqn:stab_E_h} of the discrete harmonic extension, we deduce that
$$\|E_h(\Pi_h\dot\psi_k)-E_h\dot\psi_{k,h}\|_{H^1(\Omega)}\le D_E\|\Pi_h\dot\psi_k-\dot\psi_{k,h}\|_{H^{1/2}(\Gamma)}.$$
Therefore, it remains to estimate the finite element error $\Pi_h\dot\psi_k-\dot\psi_{k,h}$ in the $H^{1/2}(\Gamma)$ norm. 
The approximation of $\dot \psi_{k}$ by $\dot \psi_{k,h}$ gives rise to two sources of error: the approximation of the source  term $\mu= \lambda \circ E$ by $\mu_h=\lambda \circ E_h$ and the approximation of the reaction-diffusion problems using a finite element method. To separate the two errors, we introduce the intermediate function $\widetilde\psi_k:=\mathcal{Q}_k^{-s}(\mathcal L)\widetilde\mu$, where $\widetilde\mu:=\lambda\circ E_h\circ \Pi_h$ and $\Pi_h$ is the Scott-Zhang projection. 
Note that proceeding as in \eqref{eqn:reg_mu} but using \eqref{eqn:stab_E_h} as well as \eqref{eq:stab_pi_h} with $t=1/2$, we deduce that $\widetilde \mu \in H^{-1/2}(\Gamma)$ with
\begin{equation} \label{eqn:bound_mu_tilde}
 \| \widetilde \mu \|_{H^{-1/2}(\Gamma)} \lesssim \| \lambda \|_{H^1(\Omega)^*}.   
\end{equation}

We write
\begin{align} \label{eqn:III_split}
  \|\Pi_h\dot \psi_k-\dot \psi_{k,h}\|_{H^{1/2}(\Gamma)} \le &\|\Pi_h(\dot \psi_k-\widetilde\psi_k)\|_{H^{1/2}(\Gamma)}+\|\Pi_h\widetilde\psi_k-\widetilde\psi_k\|_{H^{1/2}(\Gamma)} \\
  &+\|\widetilde\psi_k-\dot \psi_{k,h}\|_{H^{1/2}(\Gamma)} \nonumber
\end{align} 
and estimate each term separately.
Note that for the last term, $\widetilde\psi_k$ and $\dot \psi_{k,h}$ are defined with the same data $\widetilde \mu = \lambda \circ E_h \circ \Pi_h$ because $\mu_h(v_h) = \widetilde \mu(v_h)$ for all $v_h \in \mathbb T_h$. Moreover, Lemma~\ref{lem:boundII} ensures that $\dot \psi_k$ and $\widetilde\psi_k$ belong to $H^{1/2+\widetilde\gamma-\epsilon}(\Gamma)$. Whence, Lemma~\ref{prop:err_uk_ukh}  and \eqref{eqn:bound_mu_tilde} guarantee that
$$
\|\widetilde\psi_k-\dot \psi_{k,h}\|_{H^{1/2}(\Gamma)}\lesssim C(\epsilon) h^{\widetilde\gamma-\epsilon} \| \widetilde  \mu \|_{H^{-1/2}(\Gamma)}\lesssim C(\epsilon) h^{\widetilde\gamma-\epsilon} \| \lambda\|_{H^{1}(\Omega)^*}.
$$
The second term in \eqref{eqn:III_split} is an interpolation error and can be estimated as in \eqref{e:stb_psi} and \eqref{eqn:II_case1a} yielding
$$\|\Pi_h\widetilde\psi_k-\widetilde\psi_k\|_{H^{1/2}(\Gamma)}\lesssim C(\epsilon) h^{\widetilde\gamma-\epsilon} \| \widetilde  \mu \|_{H^{-1/2}(\Gamma)}\lesssim C(\epsilon) h^{\widetilde\gamma-\epsilon} \| \lambda\|_{H^{1}(\Omega)^*}.$$

We now focus on the first term in \eqref{eqn:III_split} and invoke the stability \eqref{eq:stab_pi_h} with $t=1/2$ to obtain
$$\|\Pi_h(\dot \psi_k-\widetilde\psi_k)\|_{H^{1/2}(\Gamma)}\lesssim \|\dot \psi_k-\widetilde\psi_k\|_{H^{1/2}(\Gamma)}.$$
Moreover, Lemma~\ref{prop_stability_Qks} with $g=\mu-\widetilde\mu$, $\beta=1/2$ and $r=\widetilde\gamma+1/2-\epsilon/2$ yields
\begin{align*}
  \|\dot \psi_k-\widetilde\psi_k\|_{H^{1/2}(\Gamma)} &\lesssim \frac{1}{\epsilon}\|\mu-\widetilde\mu\|_{H^{-r}(\Gamma)}=\frac{1}{\epsilon}\sup_{v \in H^{r}(\Gamma) \setminus \{0\}} \frac{\lambda(Ev)-\lambda(E_h(\Pi_h v))}{\| v \|_{H^{r}(\Gamma)}} \\
  &\lesssim \frac{1}{\epsilon}\| \lambda\|_{H^1(\Omega)^*} \sup_{v \in H^{r}(\Gamma) \setminus \{0\}}\frac{\|Ev-E_h(\Pi_h v)\|_{H^{1}(\Omega)}}{\| v \|_{H^{r}(\Gamma)}}.
\end{align*}
Now proceeding as in the proof of Lemma~\ref{lem:boundII} we infer that 
$$\|Ev-E_h(\Pi_h v)\|_{H^{1}(\Omega)} \lesssim C(\epsilon)h^{\min(r-1/2-\epsilon/2,\alpha_\Omega)} \| v \|_{H^r(\Gamma)}$$
and thus
$$
\|\dot \psi_k-\widetilde\psi_k\|_{H^{1/2}(\Gamma)} \lesssim C(\epsilon) h^{\min(\widetilde\gamma-\epsilon,\alpha_\Omega)} \| \lambda \|_{H^1(\Omega)^*}. 
$$

The desired estimate follows upon using the estimates for each term in the right-hand side of \eqref{eqn:III_split}.

\end{proof}

\section{Numerical experiments}
\label{sec:num_res}


In this section, we illustrate the performances of the recovery algorithm (Algorithm~\ref{alg:practical}). For the computational domain, we consider the L-shaped domain $\Omega=(-1,1)^2\setminus\{[0,1)\times(-1,0]\}$. 

We consider two functions to be recovered, namely
\begin{equation} \label{eqn:smooth}
    u_1(x,y) = e^x\cos(y)
\end{equation}
and 
\begin{equation} \label{eqn:nonsmooth}
    u_2(x,y) = \widetilde u (r,\theta) = r^{\frac{2}{3}}\sin\left(\frac{2\theta}{3}\right), \quad x=r\cos(\theta), \,\, y = r\sin(\theta),
\end{equation}
referred to as the smooth and non-smooth solutions, respectively. 
Note that although $u_1$ is smooth, its trace belongs to any $H^s(\Gamma)$, $s<3/2$ because $\Gamma$ is not $C^1$. Also $u_2$ is only in $H^{5/3}(\Omega)$ and thus $u_2|_\Gamma\in H^{7/6}(\Gamma)$. In both cases, we have $f=0$ as the functions are harmonic in $\Omega$, and thus $u_f=0$ in \eqref{e:split}. 

The measurements are Gaussian functionals mimicking point evaluations, that is of the form
\begin{equation} \label{def:lambda}
    \lambda(v;z) = \frac{1}{\sqrt{2\pi r^2}}\int_{\Omega}\exp{\left(-\frac{|x-z|^2}{2r^2}\right)}v(x){\rm d}x, \quad v\in H^1(\Omega), \quad \textrm{with }r=0.1.
\end{equation}
The centers $z$ are uniformly distributed in $\Omega$: given an integer $p$ controlling the number of points in each direction, we consider the centers
$$z_{i,j}=(-1+i\delta_p,-1+j\delta_p) \quad \text{for} \quad 1\le i,j\le p, \quad \text{where} \quad \delta_p=\frac{2}{p+1},$$
as long as $z_{i,j} \in \Omega$. Table~\ref{t:p_m} reports the number of measurements for the values of $p$ used in the numerical experiments.

\begin{table}[ht!]
\begin{center}
\begin{tabular}{c||c|c|c|c|c|c}
    $p$ & 2 & 3 & 4 & 5 & 6 & 7 \\
\hline    
    $m$ & 3 & 5 & 12 & 16 & 27 & 33 \\
\end{tabular}
\caption{Number of measurements for the values of $p$ used in the numerical experiments.}\label{t:p_m}
\end{center}
\end{table}

As mentioned earlier, the functions to be recover belong to $H^s(\Gamma)$ for $s<3/2$ (smooth case) and $s=7/6$ (rough case). In particular, both have $s>1$. 
We now briefly describe how to extend our Algorithm~\ref{alg:practical} for $s>1$.
When $s=l+\overline{s}$ with $l \geq 1$ integral and $0 \leq  \overline{s}<1$, we obtain an approximation of $\psi = \mathcal L^{-s}\mu$ by proceeding iteratively as described in Algorithm~\ref{algo:s_large}.
\begin{algorithm}
\caption{Fractional diffusion problem for $s \geq 1$}\label{algo:s_large}
\begin{algorithmic}
\Require $s=l+\overline{s}$ with $l \geq 1$ integral and $0 \leq  \overline{s}<1$
\Ensure $\dot \psi_{k,h}$
\State $\diamond$ Set $\psi^{(0)}_{h}= \mu_h$
\For{$i=1,2,\ldots,l$} 
\State $\diamond$ Compute $\psi^{(i)}_{h}=\mathcal{L}_h^{-1}\psi^{(i-1)}_{h}$
\EndFor
\If{$\overline{s}>0$}
\State $\diamond$ Compute $\dot \psi_{k,h} = \mathcal Q^{-\overline{s}}_k(\mathcal L_h) \psi^{(l)}_h$
\Else 
\State $\diamond$ $\dot \psi_{k,h}=\psi^{(l)}_h$
\EndIf
 \end{algorithmic}
\end{algorithm}

Regarding the discretization parameters, we use uniform partitions of $\Omega$ made of $3\cdot 2^{2n}$ quadrilaterals of side-length $2^{-n}$, where $n=2,\ldots,9$ is the refinement level, corresponding to a mesh size of $h=\sqrt{2} \cdot 2^{-n}$. Also, in all experiments, the sinc quadrature parameters $\texttt{M}$ and $\texttt{N}$ are taken sufficiently large not to influence the results.

The recovery functions constructed by Algorithm~\ref{alg:OR} and associated to $u_i$, $i=1,2$,  are denoted $\hat u_i$. 

\subsection{Performances of Algorithm~\ref{alg:OR} for the recovery of $u$}

We analyze the behavior of the relative recovery error 
$$e_i:=\frac{\|u_i-\widehat u_i\|_{H^1(\Omega)}}{\|u_i\|_{H^1(\Omega)}}$$
for both the smooth ($i=1$) and non-smooth ($i=2$) for different values of $s$ (regularity on the boundary), $m$ (number of measurements), and $n$ (refinement level). Note that we have $\|u_1\|_{H^1(\Omega)} \approx 2.648$ for the smooth solution and $\|u_2\|_{H^1(\Omega)} \approx 1.709$ for the non-smooth one.

In view of Theorem~\ref{t:main}, the recovery error is expected to reduce as $n$ increases until the optimal recovery threshold (Chebyshev radius $R(\dot K^s_{\boldsymbol \omega})$) is reached and as $s$ increases (at least until the function to recover belongs to $\dot K^s_{\boldsymbol \omega}$). We also expect that this threshold decreases as $m$ increases. 
This is indeed what Table~\ref{tab:recovery_smooth} (smooth case) and Table~\ref{tab:recovery_nonsmooth} (non-smooth case) indicate.
Furthermore, we observe that the recovery error increase as $m$ increase for $n$ fixed. 
This effect is already documented in \cite{BBCDDP2023} and reflects the influence of the measurements in the constant $\dot \Lambda$ in Corollary~\ref{c:main}.
We also note that the relative recovery error does not seem to be affected by the smoothness of the function except that larger $s$ can be chosen in the smooth case.

\begin{table}[htbp]
    \centering
    \begin{tabular}{|c|c|c|c|c|c|c|c|c|c|c|}
\cline{4-11}
\multicolumn{3}{c}{ } & \multicolumn{4}{|c|}{$s=1$} & \multicolumn{4}{|c|}{$s=1.45$} \\
\hline
$n$ & $\#$ cells & $\#$ dofs & $m=3$ & $m=12$ & $m=16$ & $m=33$ & $m=3$ & $m=12$ & $m=16$ & $m=33$ \\
\hline
2 & 48 & 65 & 0.71660 & 0.23765 & 0.95774 & 1264.02 & 0.63811 & 0.17785 & 1.06128 & 218.399 \\
\hline
3 & 192 & 225 & 0.70115 & 0.25950 & 0.11432 & 0.53311 & 0.62669 & 0.18755 & 0.08289 & 0.60300 \\
\hline
4 & 768 & 833 & 0.69910 & 0.26335 & 0.11677 & 0.093724 & 0.62505 & 0.18931 & 0.08058 & 0.08118 \\
\hline
5 & 3072 & 3201 & 0.69838 & 0.26429 & 0.11762 & 0.08191 & 0.62445 & 0.18965 & 0.08020 & 0.06073 \\
\hline
6 & 12288 & 12545 & 0.69813 & 0.26451 & 0.11785 & 0.08176 & 0.62423 & 0.18969 & 0.08010 & 0.06029 \\
\hline
7 & 49152 & 49665 & 0.69803 & 0.26455 & 0.11790 & 0.08181 & 0.62415 & 0.18968 & 0.08008 & 0.06033 \\
\hline
\end{tabular}
    \caption{Normalized recovery error $e_1$ for different refinement levels and different number of Gaussian measurements for $s=1$ and $s=1.45$ in the smooth solution case \eqref{eqn:smooth}.}
    \label{tab:recovery_smooth}
\end{table}

\begin{table}[htbp]
    \centering
    \begin{tabular}{|c|c|c|c|c|c|c|c|c|c|c|}
\cline{4-11}
\multicolumn{3}{c}{ } & \multicolumn{4}{|c|}{$s=0.55$} & \multicolumn{4}{|c|}{$s=7/6$} \\
\hline
$n$ & $\#$ cells & $\#$ dofs & $m=3$ & $m=12$ & $m=16$ & $m=33$ & $m=3$ & $m=12$ & $m=16$ & $m=33$ \\
\hline
2 & 48 & 65 & 0.74285 & 0.27451 & 1.38862 & 6871.97 & 0.34158 & 0.17675 & 1.49918 & 1047.31 \\
\hline
3 & 192 & 225 & 0.70773 & 0.22492 & 0.19341 & 0.94360 & 0.34456 & 0.11568 & 0.12992 & 0.97837 \\
\hline
4 & 768 & 833 & 0.69627 & 0.22385 & 0.17725 & 0.19350 & 0.34271 & 0.10777 & 0.07942 & 0.16295 \\
\hline
5 & 3072 & 3201 & 0.69128 & 0.22290 & 0.17574 & 0.13210 & 0.34190 & 0.10417 & 0.07067 & 0.06218 \\
\hline
6 & 12288 & 12545 & 0.68907 & 0.22238 & 0.17521 & 0.12761 & 0.34156 & 0.10265 & 0.06848 & 0.05098 \\
\hline
7 & 49152 & 49665 & 0.68811 & 0.22214 & 0.17494 & 0.12694 & 0.34141 & 0.10203 & 0.06769 & 0.04937 \\
\hline
\end{tabular}
    \caption{Normalized recovery error $e_2$ for different refinement levels and different number of Gaussian measurements for $s=0.55$ and $s=7/6$ in the non-smooth solution case \eqref{eqn:nonsmooth}.}
    \label{tab:recovery_nonsmooth}
\end{table}

\begin{figure}[htbp]
    \centering
    \includegraphics[width=5.1cm]{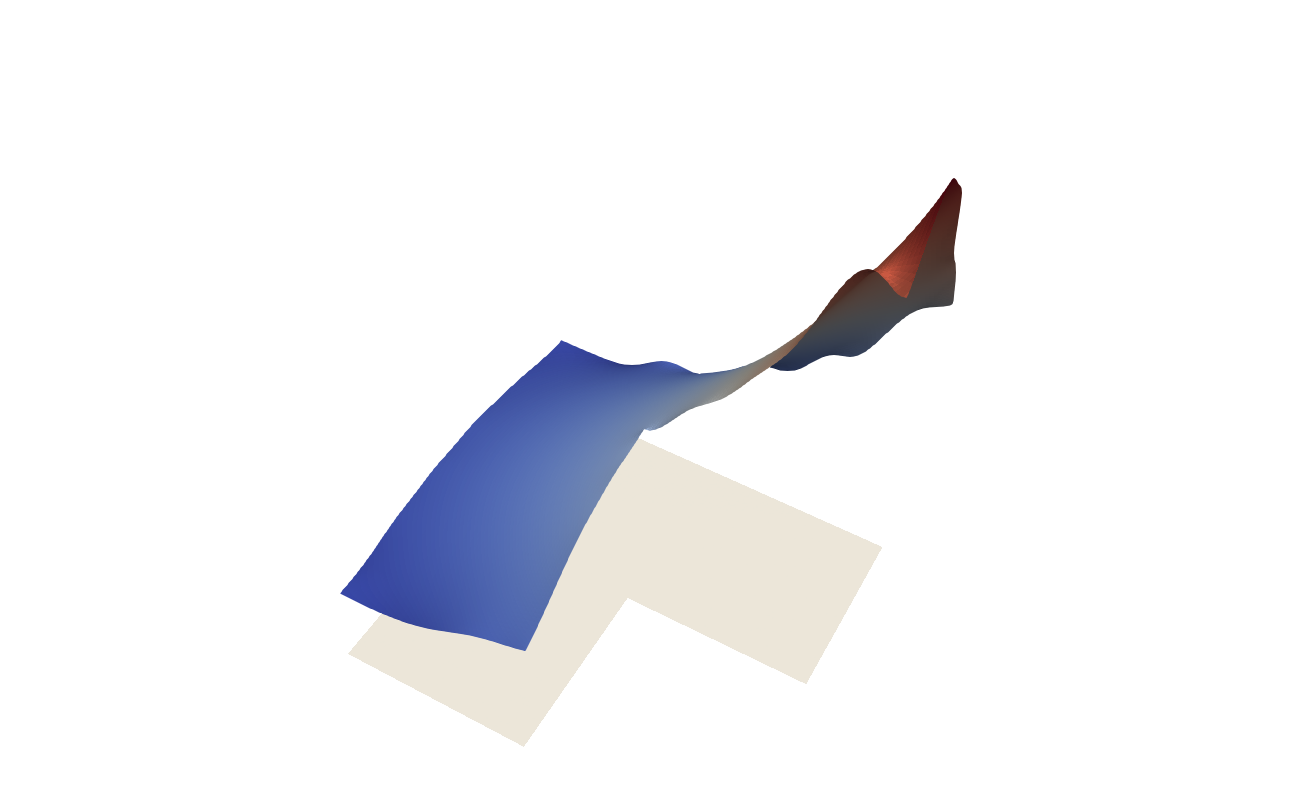}
    \hspace*{-2cm}\includegraphics[width=5.1cm]{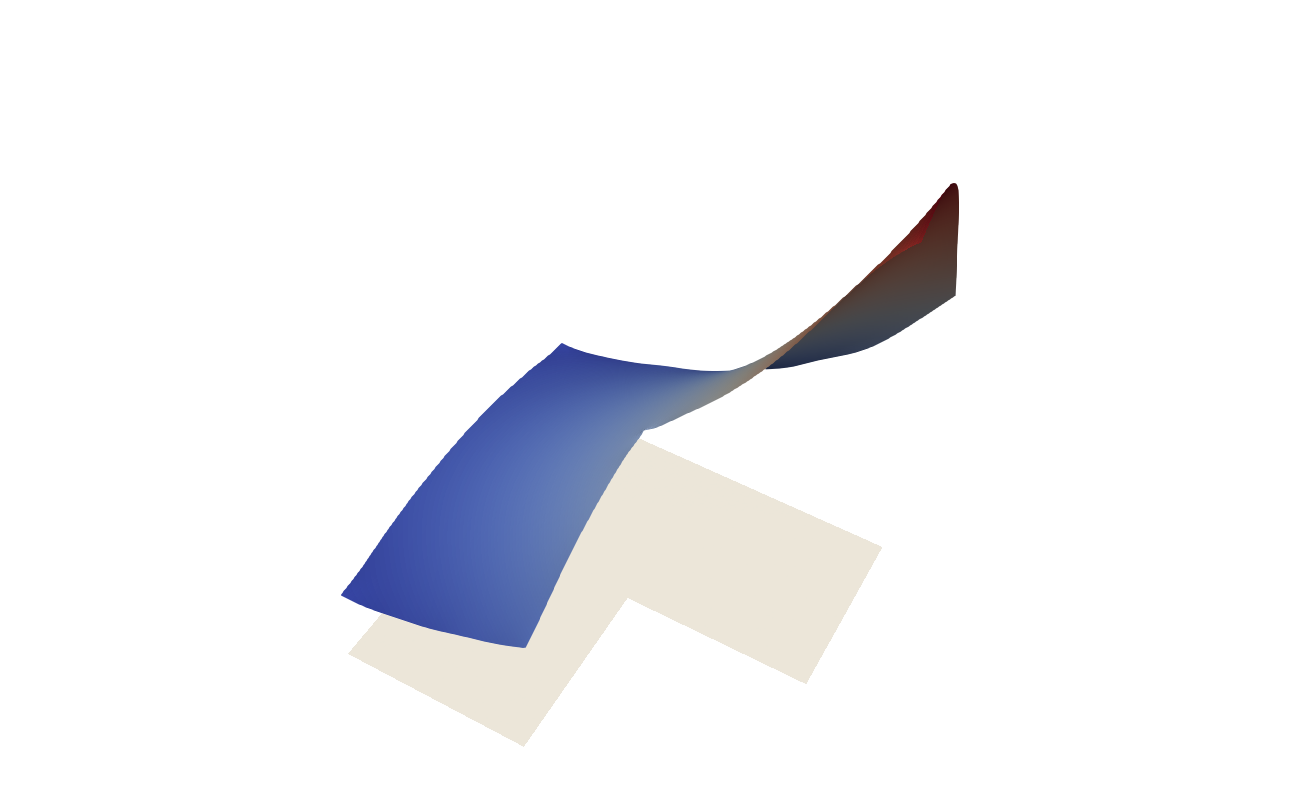}
    \hspace*{-2cm}\includegraphics[width=5.1cm]{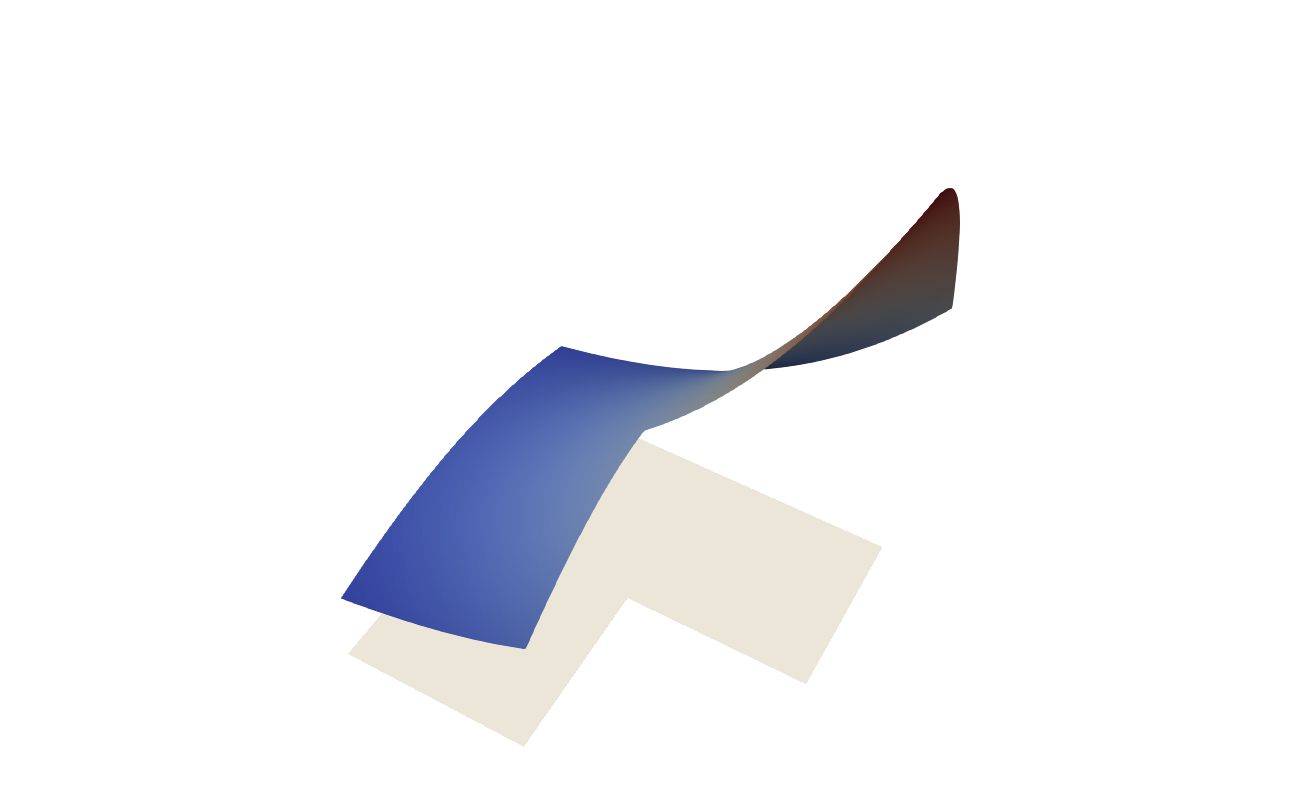}
    \caption{Recovery of the smooth solution $u_1$ in \eqref{eqn:smooth} with $m=33$ Gaussian measurements. Left: $s=0.55$ ($e_1=0.17372$), middle: $s=1.45$ ($e_1=0.06029$), right: $u_1$.}
    \label{fig:smooth_s}
\end{figure}

\begin{figure}[htbp]
    \centering
    \includegraphics[width=5.1cm]{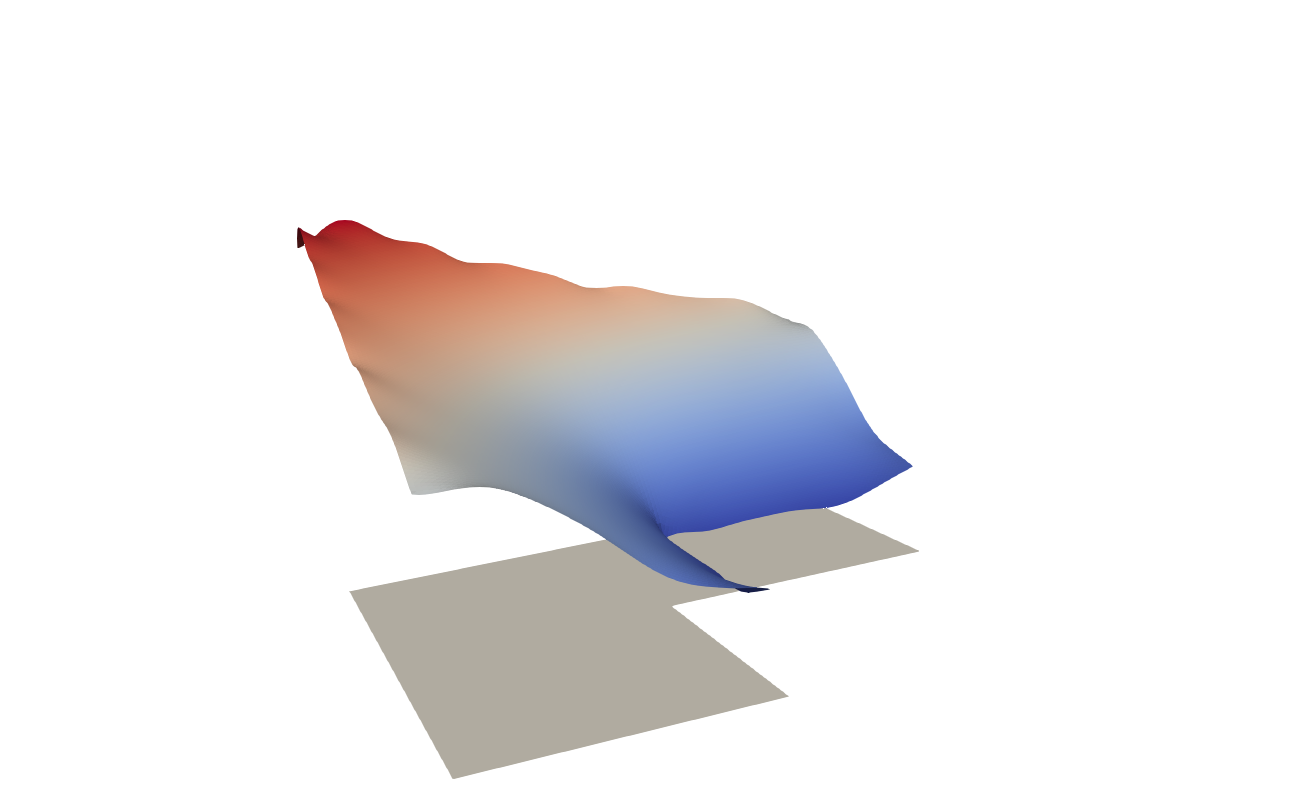}
    \hspace*{-2cm}\includegraphics[width=5.1cm]{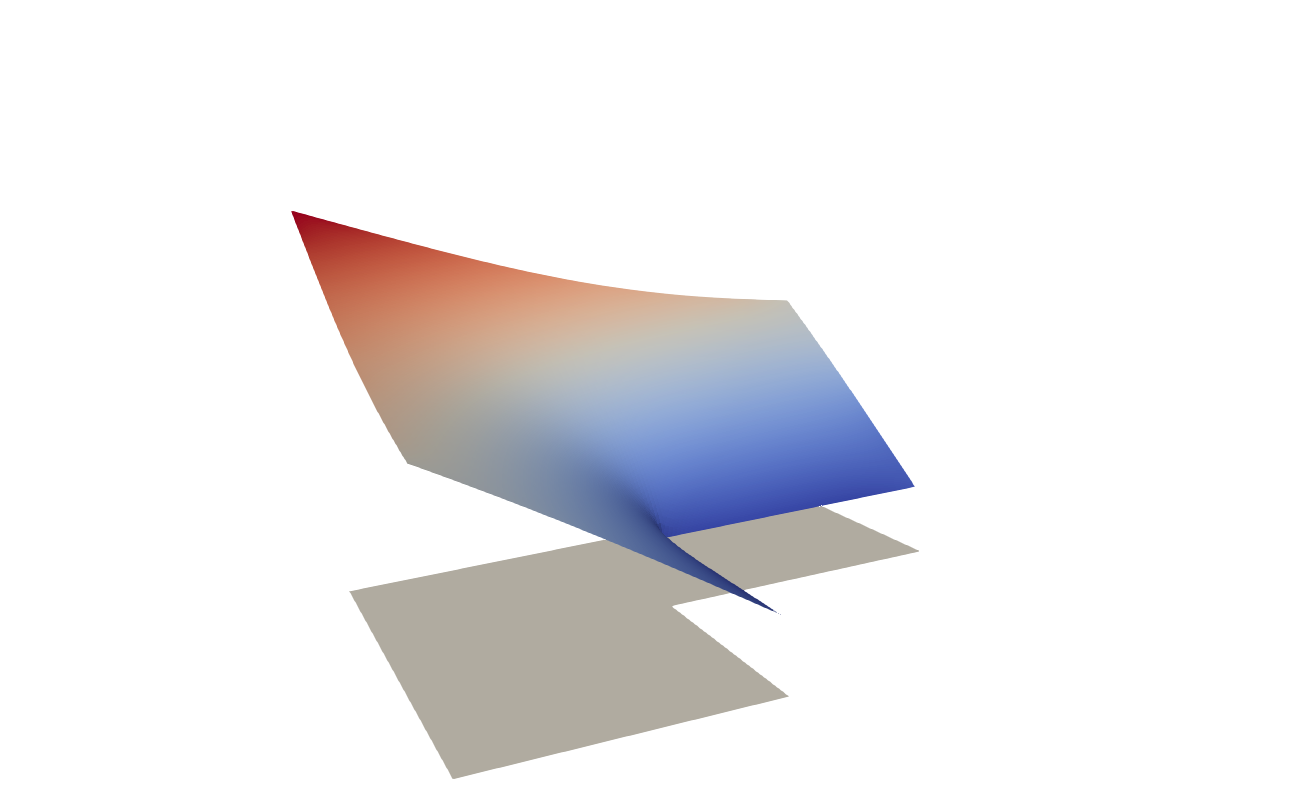}
    \hspace*{-2cm}\includegraphics[width=5.1cm]{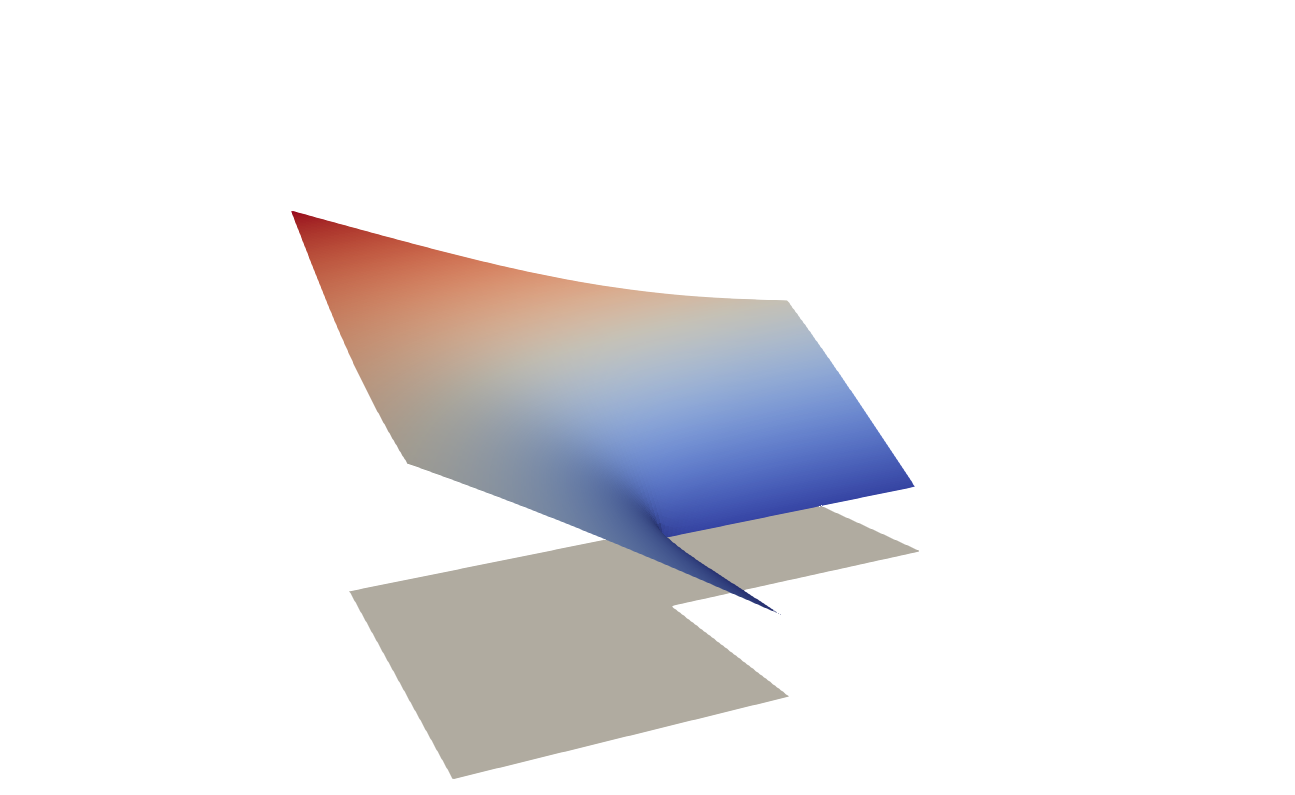}    
    \caption{Recovery of the non-smooth solution $u_2$ in \eqref{eqn:nonsmooth} with $m=33$ Gaussian measurements. Left: $s=0.55$ ($e_2=0.12761$), middle: $s=7/6$ ($e_2=0.05098$), right: $u_2$.}
    \label{fig:non_smooth_s}
\end{figure}

\begin{figure}[htbp]
    \centering
    \includegraphics[width=5.1cm]{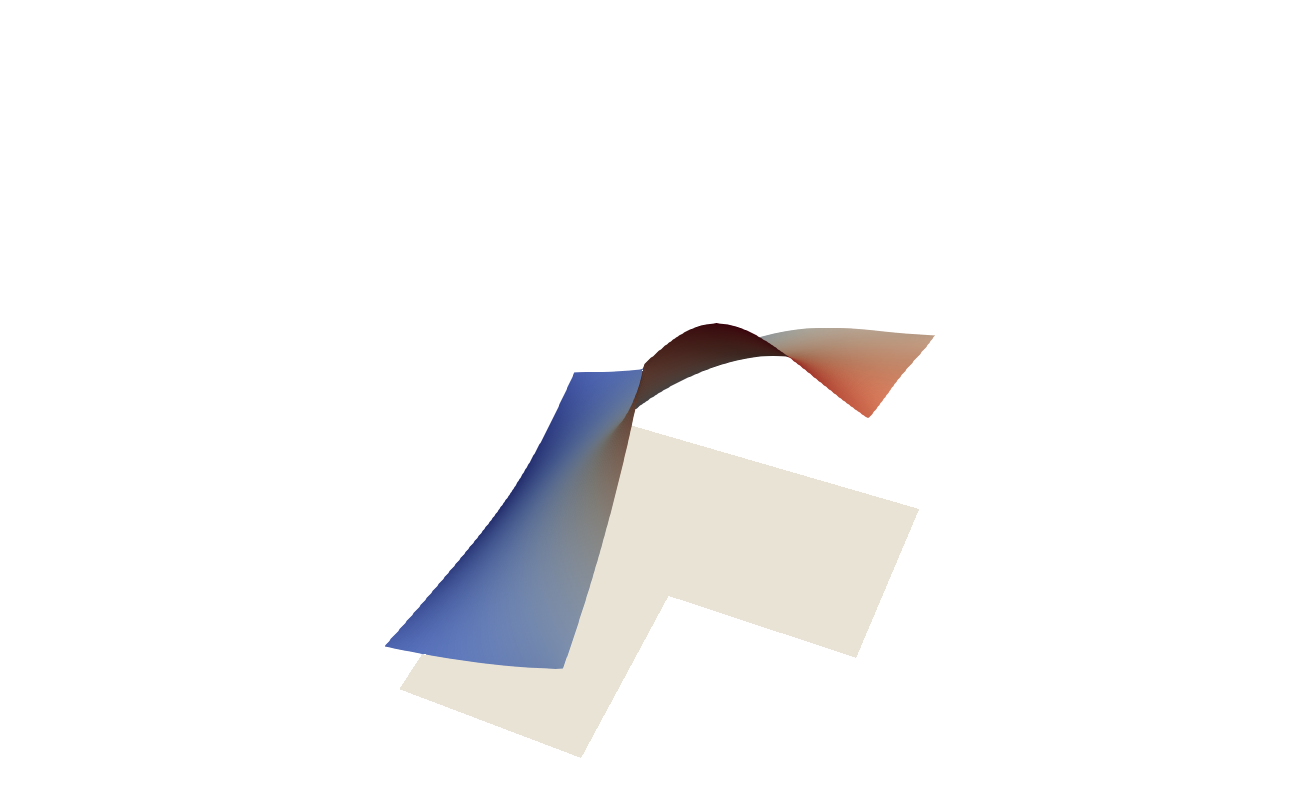}
    \hspace*{-2cm}\includegraphics[width=5.1cm]{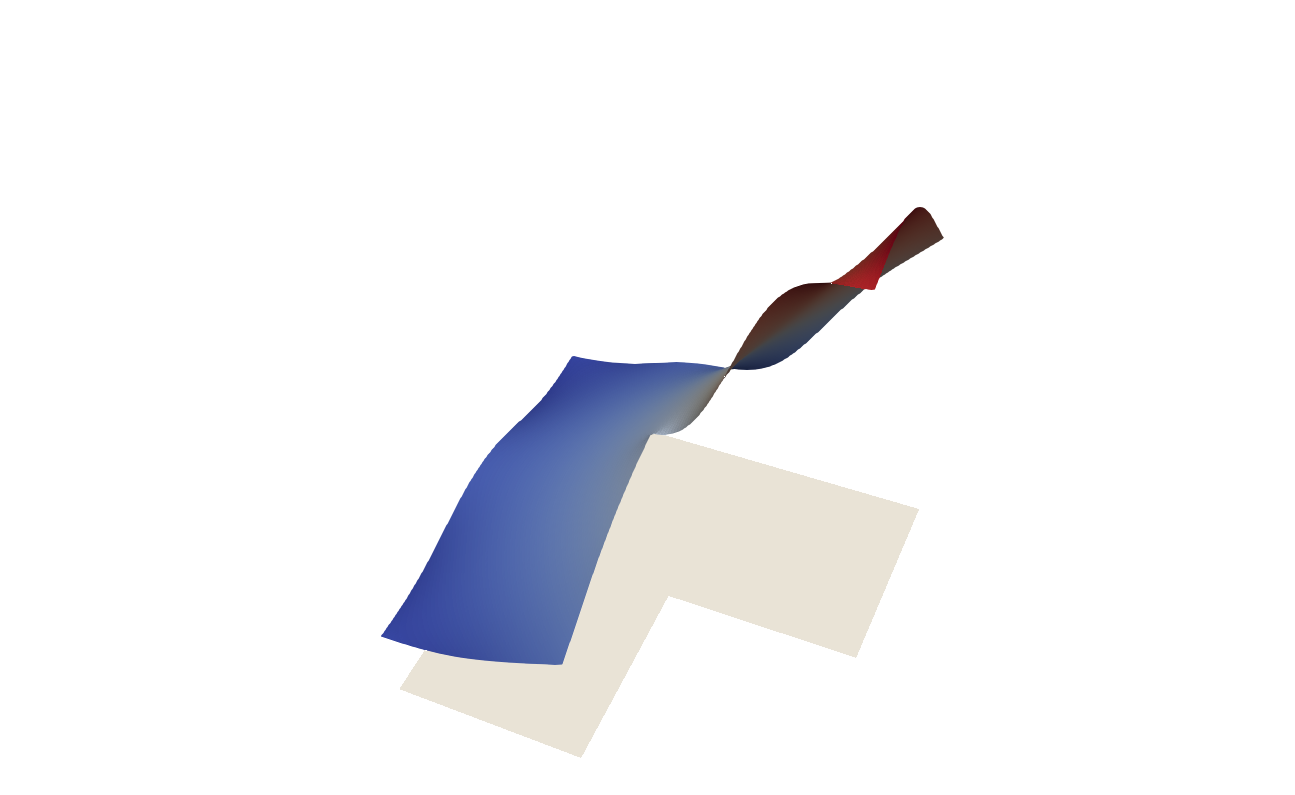}
    \hspace*{-2cm}\includegraphics[width=5.1cm]{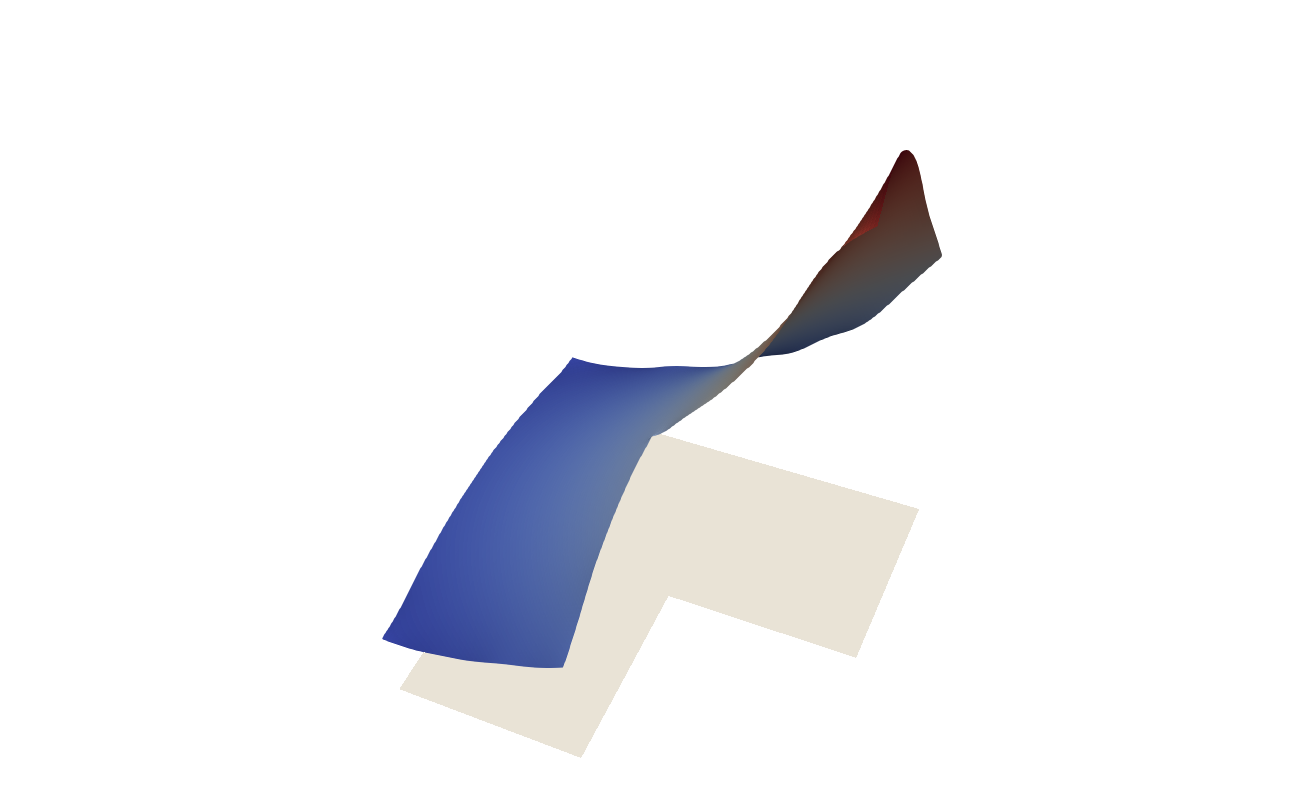}
    \caption{Recovery of the smooth solution $u_1$ in \eqref{eqn:smooth} with $s=1$ and $m$ Gaussian measurements. Left: $m=3$ ($e_1=0.69813$), middle: $m=12$ ($e_1=0.26451$), right: $m=33$ ($e_1=0.08176$).}
    \label{fig:smooth_m}
\end{figure}

\begin{figure}[htbp]
    \centering
    \includegraphics[width=5.1cm]{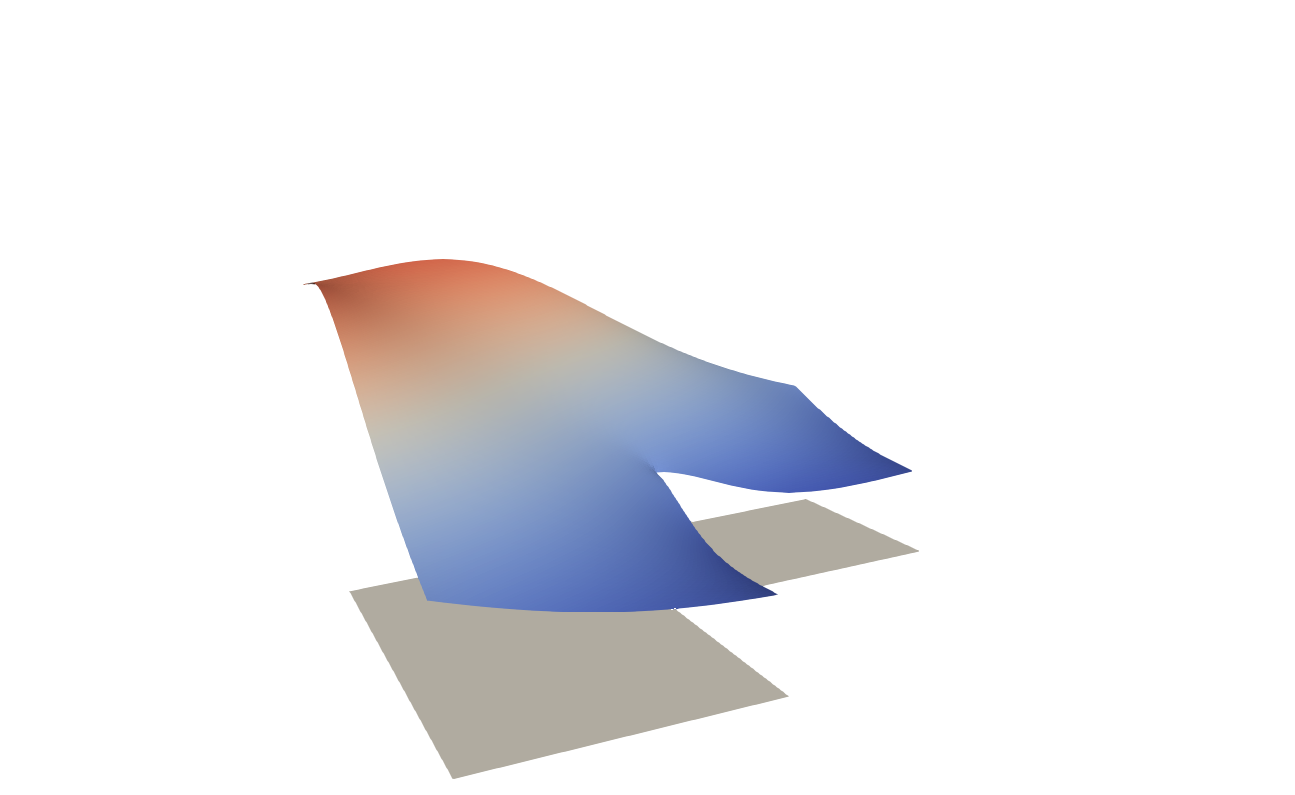}
    \hspace*{-2cm}\includegraphics[width=5.1cm]{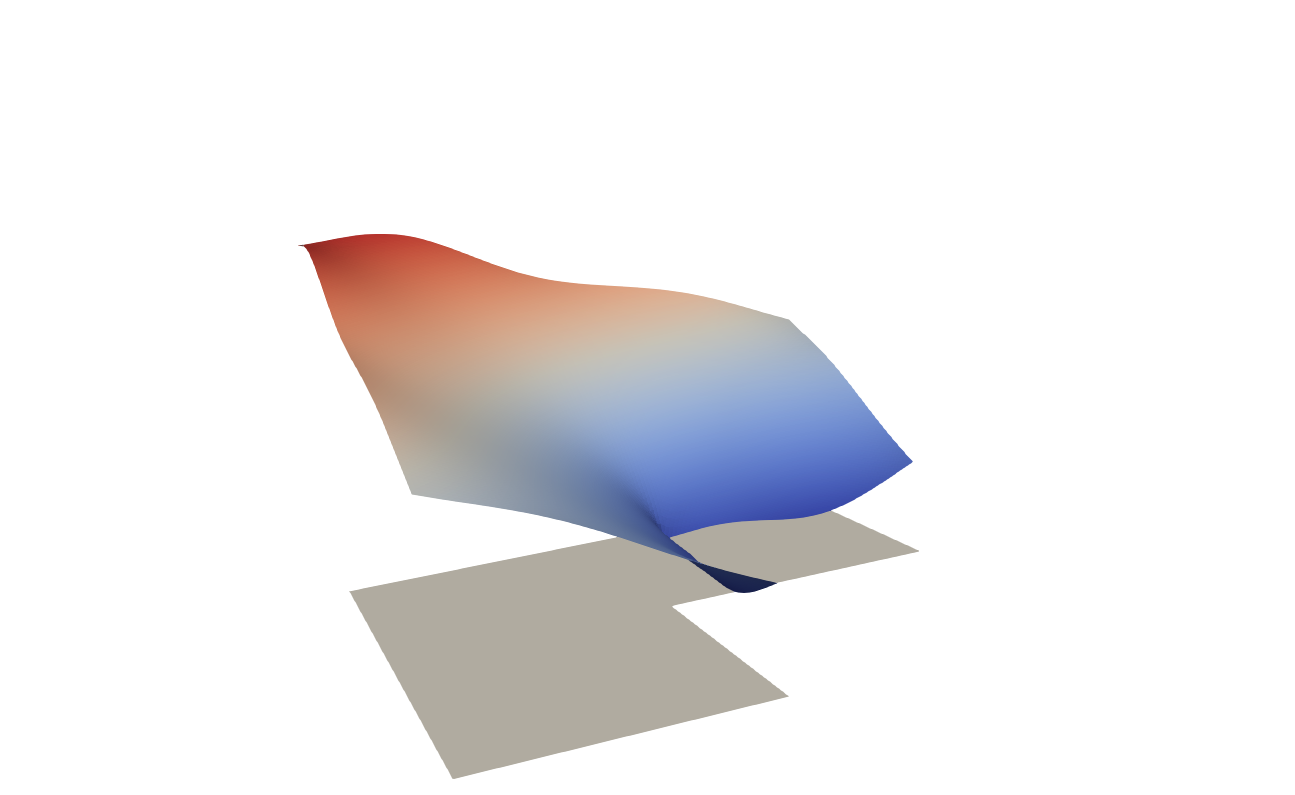}
    \hspace*{-2cm}\includegraphics[width=5.1cm]{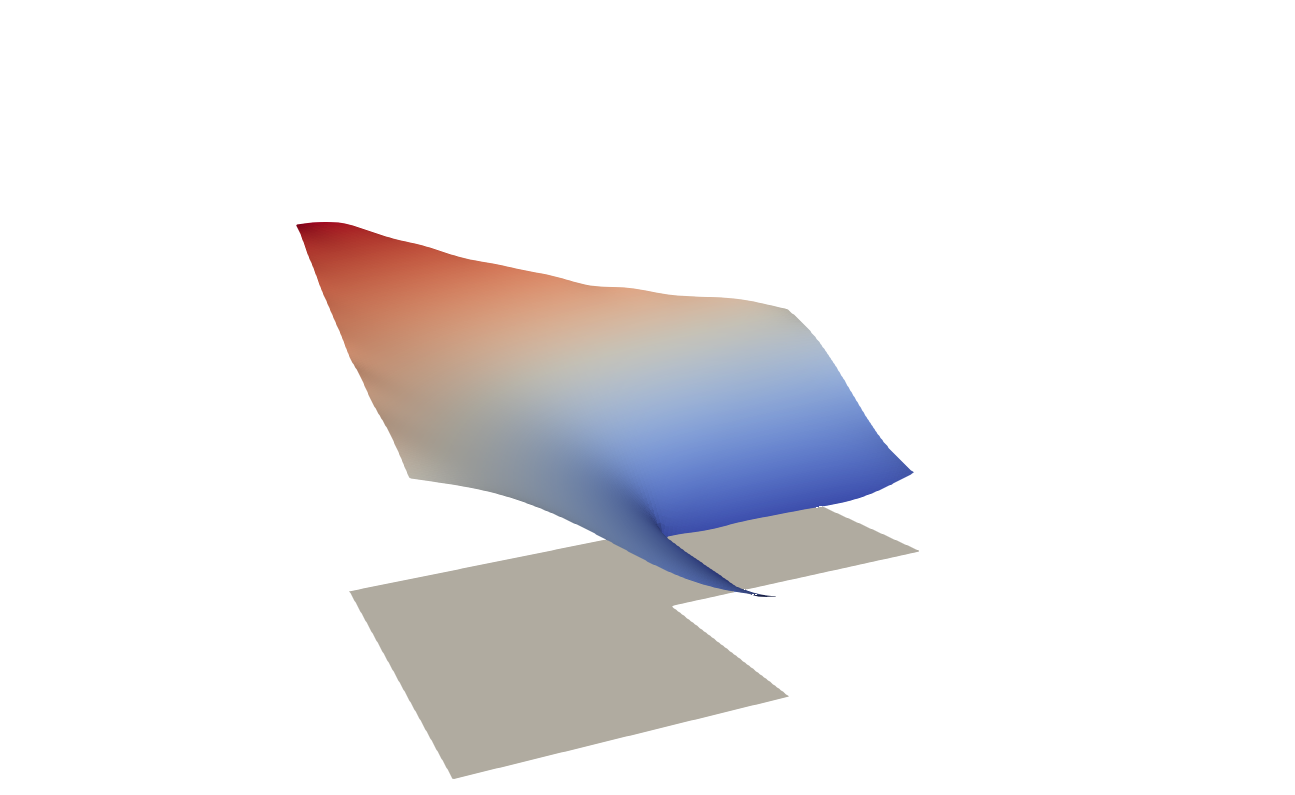}
    \caption{Recovery of the non-smooth solution $u_2$ in \eqref{eqn:nonsmooth} with $s=1$ and $m$ Gaussian measurements. Left: $m=3$ ($e_2=0.40164$), middle: $m=12$ ($e_2=0.11890$), right: $m=33$ ($e_2=0.05865$).}
    \label{fig:non_smooth_m}
\end{figure}

The influence of $s$ in the recovered functions is illustrated in Figures~\ref{fig:smooth_s} and~\ref{fig:non_smooth_s} when using $m=33$ measurements and $n=6$.
For small $s$ the recovered solution fails to approximate accurately the exact solution near the boundary as indicated by the presence of spurious oscillations. However, the latter are reduced as $s$ increases.
Instead, Figures~\ref{fig:smooth_m} and~\ref{fig:non_smooth_m} depict the recovered functions for $s=1$ and $n=6$, but for different values of the number of measurements $m$.

We end this section with an open question: what can we expect from Algorithm~\ref{alg:OR} if the $s$ chosen is larger than the regularity of the functions to be recovered? 
In Table~\ref{tab:recovery_s}, we report the recovery error for different values of $m$ and $s$ going beyond the limiting cases ($s<3/2$ for $u_1$ and $s<7/6$ for $u_2$).  
We observe that increasing $s$ beyond the value dictated by the regularity of $u_i$ continues to yield smaller recovery errors but only until a certain limiting value after which the recovery errors increase. Whether this effect is a numerical artifact or can be predicted mathematically is left open for future studies. 

\begin{table}[htbp!]
    \centering
    \begin{tabular}{|r|c|c|c|c|}
\cline{2-5}
\multicolumn{1}{c}{ } & \multicolumn{2}{|c|}{ $u_1$ (smooth)} & \multicolumn{2}{|c|}{$u_2$ (non-smooth)} \\
\hline
\multicolumn{1}{|c|}{$s$} & $m=3$ & $m=33$ & $m=3$ & $m=33$ \\
\hline
0.55 & 0.85958 & 0.17372 & 0.68908 & 0.12761 \\
0.66 & 0.80998 & 0.13560 & 0.59458 & 0.09864 \\
0.75 & 0.77457 & 0.11460 & 0.53124 & 0.08277 \\
1.00 & 0.69813 & 0.08176 & 0.40164 & 0.05865 \\
7/6  &0.66269  &0.07059 & 0.34156  &0.05098 \\
1.25 & 0.64885 & 0.06671 & 0.31710 & 0.04845 \\
1.50 & 0.61967 & 0.05915 & 0.26038 & 0.04382 \\
1.75 & 0.60393 & 0.05512 & 0.22215 & 0.04170 \\
1.95 & 0.59755 & 0.05325 & 0.20102 & 0.04093 \\
2.00 & 0.59656 & 0.05290 & 0.19677 & 0.04082 \\
5.00 & 0.61386 & 0.08286 & 0.15607 & 0.04479 \\
10.00&  0.61994&  0.16415& 0.15670 & 0.11683 \\
20.00&  0.62008&  2.96580& 0.15672 & 0.95191\\
\hline
\end{tabular}
    \caption{Normalized recovery error $e_i$ for different values of $s$ and number of Gaussian measurements for both the smooth ($i=1$) and non-smooth ($i=2$) cases. We use $n=6$ for the refinement level of the mesh.}
    \label{tab:recovery_s}
\end{table}

\subsection{Performances of Algorithm~\ref{alg:practical} for the approximation of the Riesz representers}

\begin{figure}[htbp]
    \centering
    \includegraphics[width=5.8cm]{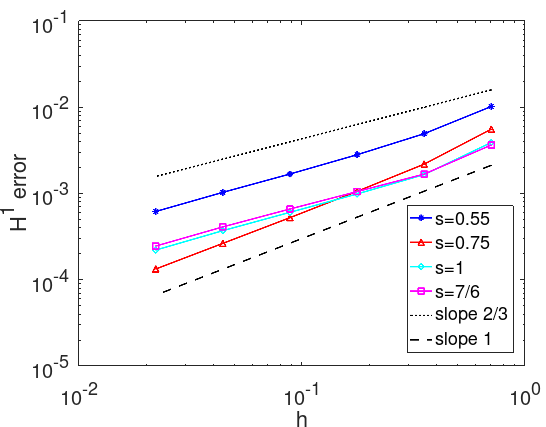}
    \includegraphics[width=5.8cm]{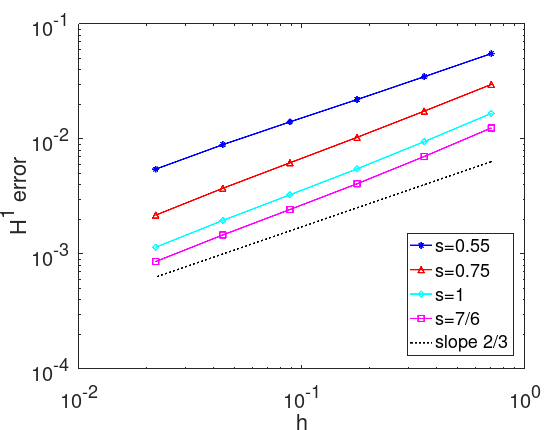}    
    \caption{Convergence of the Riesz representer for $n=2,...,6$. The $H^1(\Omega)$ error is computed using a reference discrete solution ($n=9$). Left: Gaussian measurement centered at $z=(-0.5,0.5)$; right: Gaussian measurement centered at $z=(-0.1,0.1)$.}
    \label{fig:Riesz_convergence}
\end{figure}

We end the numerical section with a discussion on the approximation of the Riesz representers using Algorithm~\ref{alg:practical} for two different choices of centers $z$ in the Gaussian measurements \eqref{def:lambda}. In absence of an exact solution, we compute the $H^1(\Omega)$ error between the output of Algorithm~\ref{alg:practical} and a reference discrete solution obtained using $n=9$ refinements ($786432$ cells). The number of sinc quadrature points is large enough not to influence the approximations. 
The $H^1(\Omega)$ error for different values of the mesh size $h$ is provided in Figure~\ref{fig:Riesz_convergence}. In most cases, the observed convergence rate is approximately $0.75$, which is marginally better than $2/3$. Note that for the L-shaped domain $\alpha_\Omega = 2/3-\epsilon$ for any $\epsilon>0$ and seems to dictate the convergence rate guaranteed by Theorem~\ref{t:main}. The case $z=(-0.5,0.5)$ and $s=0.75$ is an outlier for which a rate closer to $1$ is observed.

\begin{acknowledgement}

AB is partially supported by the NSF Grant DMS-2409807. DG acknowledges the support provided by the Natural Science and Engineering Research Council (NSERC, grant RGPIN-2021-04311).
\end{acknowledgement}
%

\section*{Appendix}

This appendix contains several lemmas that cannot be directly found in the literature because of the low regularity of the source terms (measurements) as we all the boundary $\Gamma$ considered in this work. However, the proofs presented in this appendix follows those provided in \cite{BL2022,BGL2024}. Those that are only marginally different are omitted. 

Given $g\in H^{-r}(\Gamma)$, $r<1$, let $w = \mathcal Tg \in H^1(\Gamma)$ be given by \eqref{eq:def_op_L}, i.e., $w$ solves
\begin{equation} \label{eqn:w}
a(w,v):=\langle w,v\rangle_{H^1(\Gamma)}=\int_{\Gamma}wv+\int_{\Gamma}\nabla_{\Gamma}w\cdot\nabla_{\Gamma}v = g(v), \qquad \forall\, v\in H^1(\Gamma).
\end{equation}

Let $w_h\in\mathbb{T}_h$ be the finite element approximation to $w$, i.e., the solution to
\begin{equation} \label{eqn:w_h}
a(w_h,v_h) = g(v_h), \qquad \forall\, v\in \mathbb{T}_h.    
\end{equation}
We start with estimates on the discrepancy $w-w_h$ including the two specific cases consistent with the estimates in \cite{BL2022,BGL2024} under Assumption~\ref{a:pickup}
\begin{itemize}
    \item $\|w-w_h\|_{H^{\beta}(\Gamma)} \in \mathcal{O}(h^{2-r-\beta})$ for  $\beta \in [0,1]$ and $\alpha_{\Gamma}=1$;
    \item $\|w-w_h\|_{H^{\beta}(\Gamma)} \in \mathcal{O}(h^{\alpha_{\Gamma}+\min\{\alpha_{\Gamma},1-\beta\}})$ for  $\beta \in [0,1]$ and $r=0$.
\end{itemize}

\begin{lemma} \label{lemma:w_w_h}
    Let $\alpha_\Gamma$ be as in Assumption~\ref{a:pickup}. Let $g\in H^{-r}(\Gamma)$ for $r<1$, and let $w$ and $w_h$ be given by \eqref{eqn:w} and \eqref{eqn:w_h}, respectively. Then for any $\beta\in[0,1]$ we have
    $$\|w-w_h\|_{H^{\beta}(\Gamma)}\lesssim h^{\gamma+\overline{\gamma}}\|g\|_{H^{-r}(\Gamma)},$$
    where $\gamma:=\min\{\alpha_{\Gamma},1-r\}$ and $\overline{\gamma}:=\min\{\alpha_{\Gamma},1-\beta\}$.
\end{lemma}

\begin{proof}
We first prove the result for $\beta=1$. Using the Galerkin orthogonality property
$$a(w-w_h,v_h)=0, \qquad \forall\, v_h\in\mathbb{T}_h,$$
we have for $v_h\in\mathbb{T}_h$
\begin{align*}
    \|w-w_h\|_{H^1(\Gamma)}^2 &= a(w-w_h,w-w_h) = a(w-w_h,w-v_h)  \\
    &\le \|w-w_h\|_{H^1(\Gamma)}\|w-v_h\|_{H^1(\Gamma)}
\end{align*}
and thus
$$\|w-w_h\|_{H^1(\Gamma)}\le \inf_{v_h\in\mathbb{T}_h}\|w-v_h\|_{H^1(\Gamma)}.$$
Now since $w\in H^{1+\gamma}(\Gamma)$, a standard finite element approximation estimate on $\Gamma$, see for instance \cite{BDN2020}, yields
$$\inf_{v_h\in\mathbb{T}_h}\|w-v_h\|_{H^1(\Gamma)}\lesssim h^{\gamma}\|w\|_{H^{1+\gamma}(\Gamma)}$$
and because $\gamma \leq \alpha_\Gamma$, Assumption~\ref{a:pickup} implies that
$$
\|w\|_{H^{1+\gamma}(\Gamma)}\lesssim \|g\|_{H^{\gamma-1}(\Gamma)}.
$$
Whence, we obtain
$$
\|w-w_h\|_{H^1(\Gamma)} \lesssim h^{\gamma}\|g\|_{H^{\gamma-1}(\Gamma)}\lesssim h^{\gamma}\|g\|_{H^{-r}(\Gamma)},
$$
and the case $\beta=1$ is proved. 

We deduce the general case $\beta\in (0,1)$ from the case $\beta=1$ using a duality argument. Let $z\in H^1(\Gamma)$ be the solution to
$$a(v,z)=\langle w-w_h,v\rangle_{\dot H^{\beta}(\Gamma)}, \qquad \forall\, v\in H^1(\Gamma).$$
Notice that $w-w_h \in H^\beta(\Gamma)$ and recall that $\langle w-w_h,v\rangle_{\dot H^{\beta}(\Gamma)}=\langle\mathcal L^{\beta}(w-w_h),v\rangle$.
Whence, $z = \mathcal Tg$ where $g = \mathcal L^{\beta}(w-w_h) \in H^{-\beta}(\Gamma)$. Thanks to Assumption~\ref{a:pickup} and the definition of $\overline{\gamma}$, we find that $z\in H^{1+\overline{\gamma}}(\Gamma)$. Choosing $v=w-w_h$ and using the Galerkin orthogonality property, we deduce that for any $v_h\in\mathbb{T}_h$ there holds
$$
\|w-w_h\|_{\dot H^{\beta}(\Gamma)}^2 = a(w-w_h,z)=a(w-w_h,z-v_h)\le \|w-w_h\|_{H^1(\Gamma)}\|z-v_h\|_{H^1(\Gamma)}. 
$$
Invoking again a standard finite element approximation estimate on $\Gamma$, we find that
$$
\inf_{v_h\in\mathbb{T}_h}\|z-v_h\|_{H^1(\Gamma)} \lesssim h^{\overline{\gamma}}\|z\|_{H^{1+\overline{\gamma}}(\Gamma)} \lesssim 
 h^{\overline{\gamma}}\|w-w_h\|_{\dot H^t(\Gamma)}\lesssim h^{\overline{\gamma}}\|w-w_h\|_{\dot H^{\beta}(\Gamma)},    
$$
for
$$t:=-1+\overline{\gamma}+2\beta \le \beta.$$
This together with the estimate previously obtained for $\beta=1$ implies
$$\|w-w_h\|_{\dot H^{\beta}}\lesssim h^{\overline{\gamma}}\|w-w_h\|_{H^1(\Gamma)}\lesssim h^{\gamma+\overline{\gamma}}\|g\|_{H^{-r}(\Gamma)}.$$
The case $\beta \in [0,1)$ is complete in view of the equivalence $\dot H^\beta(\Gamma) \cong H^\beta(\Gamma)$, see \eqref{e:equiv_Hs}.
\end{proof}



Given $\widetilde{\mu} \in  H^{-1/2}(\Gamma)$, we now derive an upper bound for the error between 
$$
\widetilde{ \psi}_k:=\mathcal{Q}_k^{-s}(\mathcal{L})\widetilde{\mu} = \frac{k\sin(\pi s)}{\pi}\sum_{l=-\texttt{M}}^{\texttt{N}}\eta_l^{1-s}\underbrace{(\eta_lI+\mathcal{L})^{-1}\widetilde{\mu}}_{=:w(\eta_l)}
$$ 
and its finite element approximation 
\begin{equation} \label{eqn:u_sinc_FE}
    \dot{\psi}_{k,h}:=\mathcal{Q}_k^{-s}(\mathcal{L}_h) \widetilde \pi_h \widetilde{\mu} = \frac{k\sin(\pi s)}{\pi}\sum_{l=-\texttt{M}}^{\texttt{N}}\eta_l^{1-s}\underbrace{(\eta_lI+\mathcal{L}_h)^{-1}\widetilde{\pi}_h \widetilde{\mu}}_{=:w_h(\eta_l)}.
\end{equation}
Here $\widetilde \pi_h:H^{-1}(\Gamma)\rightarrow\mathbb{T}_h$ is a projection operator defined for $g\in H^{-1}(\Gamma)$ by the relations
$$
\langle\widetilde \pi_h g,v_h\rangle_{L^2(\Gamma)} = \langle g,v_h\rangle,  \qquad \forall\, v_h\in\mathbb{T}_h;
$$
recall that we have identified $L^2(\Gamma)$ with its dual. 
In particular, if $g\in L^2(\Gamma)$ then $\widetilde \pi_h g=\pi_h g$ with $\pi_h:L^2(\Gamma)\rightarrow\mathbb{T}_h$ the $L^2(\Gamma)$ orthogonal projection operator. 

To estimate the error $\|\widetilde \psi_k-\dot \psi_{k,h}\|_{H^{1/2}(\Gamma)}$, we state a result slightly more general than Lemma 4.8 in \cite{BL2022}. We omit the proof because it follows the one of Lemma~4.8 in \cite{BL2022}.

\begin{lemma} \label{lemma:w_mu_FEM}
Let $\alpha_\Gamma$ be as in Assumption~\ref{a:pickup} and let $g\in H^{-r}(\Gamma)$ with $r<1$. Moreover, for $\eta\in(0,\infty)$, let $w(\eta)=(\eta I+\mathcal{L})^{-1}g$ and $w_h(\eta)=(\eta I+\mathcal{L}_h)^{-1}\widetilde \pi_h g$. Then for any $\delta\in[r,r+1]$ with $\delta<1$ and any $\tau,\beta\in[0,1]$ with $\tau+\beta\in[0,1]$ we have
\begin{equation} \label{eqn:w_w_h_mu}
    \|\pi_h w(\eta)-w_h(\eta)\|_{H^{\beta}(\Gamma)} \lesssim \eta^{\frac{r-\delta-\tau}{2}}h^{\min\{\sigma_1,\sigma_2\}}\|g\|_{H^{-r}(\Gamma)},
\end{equation}
where
$$\sigma_1:=\min\{1-\delta,\alpha_{\Gamma}\}+1-\tau-\beta \quad \text{and} \quad \sigma_2:=\min\{1-\delta,\alpha_{\Gamma}\}+\min\{1-\tau-\beta,\alpha_{\Gamma}\}.$$
\end{lemma}

In the specific case where $r=\beta=1/2$, the estimate in Lemma~\ref{lemma:w_mu_FEM} reads: for any $\delta\in[1/2,1)$ and any $\tau\in[0,1/2]$ we have
\begin{equation} \label{eqn:w_w_h_mu12}
    \|\pi_h w(\eta)-w_h(\eta)\|_{H^{1/2}(\Gamma)} \lesssim \eta^{\frac{1/2-\delta-\tau}{2}}h^{\min\{\sigma_1,\sigma_2\}}\|g\|_{H^{-1/2}(\Gamma)},
\end{equation}
where
\begin{equation} \label{eqn:sigmas}
\sigma_1=\min\{1-\delta,\alpha_{\Gamma}\}+1/2-\tau \quad \text{and} \quad \sigma_2=\min\{1-\delta,\alpha_{\Gamma}\}+\min\{1/2-\tau,\alpha_{\Gamma}\}.
\end{equation}

We are now in position to derive an estimate for $\|\widetilde \psi_k-\dot \psi_{k,h}\|_{H^{1/2}(\Gamma)}$.
\begin{lemma} \label{prop:err_uk_ukh}
Let $\alpha_\Gamma$ be as in Assumption~\ref{a:pickup}, $s\in(1/2,1)$, and
$$\widetilde\gamma:=\min\{2s-1,\alpha_{\Gamma}+1/2,2\alpha_{\Gamma}\}.$$
Then for $\widetilde \mu \in H^{-1/2}(\Gamma)$ and for $0<\epsilon<\widetilde\gamma$, there is a constant $C(\epsilon)$ such that $C(\epsilon)\to \infty$ as $\epsilon \to 0^+$ and     
    $$\|\widetilde \psi_k-\dot \psi_{k,h}\|_{H^{1/2}(\Gamma)}\leq C(\epsilon) h^{\widetilde \gamma-\epsilon}\|\widetilde \mu \|_{H^{-1/2}(\Gamma)}.$$
\end{lemma}

\begin{proof}
Without loss of generality, we assume that $h\le 1$. Let $\alpha^*:=(\widetilde\gamma-\epsilon)/2$. We start by using triangle's inequality to decompose the error as
\begin{equation}\label{e:des}
\|\widetilde \psi_k-\dot \psi_{k,h}\|_{H^{1/2}(\Gamma)}\le \|\widetilde \psi_k-\pi_h \widetilde \psi_k\|_{H^{1/2}(\Gamma)} + \|\pi_h \widetilde \psi_k-\dot \psi_{k,h}\|_{H^{1/2}(\Gamma)}.
\end{equation}
The first term is the error due to interpolation. Proceeding as in the proof of Lemma~\ref{lem:boundII} (but using $\pi_h$ instead of $\Pi_h$), we obtain
\begin{equation}\label{e:first}
\|\widetilde \psi_k-\pi_h \widetilde \psi_k\|_{H^{1/2}(\Gamma)}\lesssim C(\epsilon)h^{\widetilde\gamma-\epsilon}\|g\|_{H^{-1/2}(\Gamma)}.
\end{equation}
For the second term, we first observe that
$$\|\pi_h \widetilde \psi_k-\dot \psi_{k,h}\|_{H^{1/2}(\Gamma)}\le\frac{\sin(\pi s)}{\pi}k\sum_{l=-\texttt{M}}^{\texttt{N}}\eta_l^{1-s}\|\pi_h w(\eta_l)-w_h(\eta_l)\|_{H^{1/2}(\Gamma)}.$$
To simplify the notation, for any $\eta\in(0,\infty)$ we write $e(\eta):=\pi_h w(\eta)-w_h(\eta)$. The rest of the proof then consists in splitting the sum into three parts, namely consider
$$\eta_l\in(0,1), \quad \eta_l\in[1,h^{-\frac{2\alpha^*}{s-1/2}}), \quad \eta_l\in[h^{-\frac{2\alpha^*}{s-1/2}},\infty),$$
and use different estimates for $\|e(\eta_l)\|_{H^{1/2}(\Gamma)}$.
\begin{itemize}
    \item[\boxed{1}] Case $\eta_l\in(0,1)$. We use \eqref{eqn:w_w_h_mu12} with $\delta=1/2$ and $\tau=0$. Since
    $$\sigma_1=\min\{1/2,\alpha_{\Gamma}\}+1/2\le \widetilde\gamma \quad \text{and} \quad \sigma_2=2\min\{1/2,\alpha_{\Gamma}\}\lesssim\widetilde\gamma$$
    we obtain
    $$\|e(\eta_l)\|_{H^{1/2}(\Gamma)}\lesssim h^{\widetilde\gamma}\|g\|_{H^{-1/2}(\Gamma)}\lesssim h^{2\alpha^*}\|g\|_{H^{-1/2}(\Gamma)},$$
    and thus
    \begin{align*}
        \frac{\sin(\pi s)}{\pi}k\sum_{0<\eta_l<1}\eta_{l}^{1-s}\|e(\eta_l)\|_{H^{1/2}(\Gamma)} &\lesssim k\sum_{0<\eta_l<1}\eta_{l}^{1-s}h^{2\alpha^*}\|g\|_{H^{-1/2}(\Gamma)} \\
        &\lesssim h^{2\alpha^*}\|g\|_{H^{-1/2}(\Gamma)}.
    \end{align*}
    \item[\boxed{2}] Case $\eta_l\in[h^{-\frac{2\alpha^*}{s-1/2}},\infty)$. Using triangle's inequality and the stability property \eqref{eq:stab_pi_h} we have
    $$\|e(\eta_l)\|_{H^{1/2}(\Gamma)}\le \|w(\eta_l)\|_{H^{1/2}(\Gamma)}+\|w_h(\eta_l)\|_{H^{1/2}(\Gamma)}.$$
    The first term can then be bounded proceeding as in the proof of Lemma~\ref{prop_stability_Qks} below, namely use \eqref{eqn:inner_w_y} and \eqref{eqn:Young} with $r=\beta=1/2$ (with $y$ replaced by $\eta_l=e^{y_l}$) together with the equivalence \eqref{e:equiv_Hs}. The bound for the second term is similar upon replacing the role of $\mathcal L$ by $\mathcal L_h$ \cite{bonito2017numerical,BL2022}.
    Whence, we infer that
$$
       \|e(\eta_l)\|_{H^{1/2}(\Gamma)}\lesssim \eta_l^{-1/2}\|g\|_{H^{-1/2}(\Gamma)}. 
$$
    Consequently, we find that
    \begin{align*}
        \frac{\sin(\pi s)}{\pi}k\sum_{\eta_l\ge h^{-\frac{2\alpha^*}{s-1/2}}}\eta_{l}^{1-s}\|e(\eta_l)\|_{H^{\beta}(\Gamma)} &\lesssim k\sum_{\eta_l\ge h^{-\frac{2\alpha^*}{s-1/2}}}\eta_{l}^{1/2-s}\|g\|_{H^{-1/2}(\Gamma)} \\
        &\lesssim h^{2\alpha^*}\|g\|_{H^{-1/2}(\Gamma)};
    \end{align*}
    where the last inequality is derived as in the proof of Theorem 4.2 in \cite{BL2022}. \item[\boxed{3}] Case $\eta_l\in[1,h^{-\frac{2\alpha^*}{s-1/2}})$. We decompose this case in three subcases and invoke \eqref{eqn:w_w_h_mu12} for different values of $\delta\in[1/2,1)$ and $\tau\in[0,1/2]$. In each case, we show that
     \begin{equation}\label{e:subcase}
     \frac{\sin(\pi s)}{\pi}k\sum_{1\le \eta_l<h^{-\frac{2\alpha^*}{s-1/2}}}\eta_{l}^{1-s}\|e(\eta_l)\|_{H^{1/2}(\Gamma)} \lesssim \frac{1}{\epsilon}h^{2\alpha^*}\|g\|_{H^{-1/2}(\Gamma)}.
     \end{equation}
    \begin{enumerate}[i)]
        \item Case $\alpha_{\Gamma}\ge 1/2$. 
        Since $1-\delta\in[0,1/2)$ and $1/2-\tau \in [0,1/2]$ we have
        $\sigma_1=\sigma_2=3/2-\delta-\tau$. We then choose $\delta$ and $\tau$ so that
        $$\delta+\tau=3/2+1-2s+\epsilon\in[1/2,3/2)$$
        in which case \eqref{eqn:w_w_h_mu12} reads
        \begin{equation*} \label{eqn:err_w_case1}
        \|e(\eta_l)\|_{H^{1/2}(\Gamma)} \lesssim \eta_l^{s-1 -\epsilon/2}h^{2s-1-\epsilon}\|g\|_{H^{-1/2}(\Gamma)}\le \eta_l^{s-1 -\epsilon/2}h^{2\alpha^*}\|g\|_{H^{-1/2}(\Gamma)}.
        \end{equation*}
        Therefore, we get
        \begin{align*}
        \frac{\sin(\pi s)}{\pi}k\sum_{1\le \eta_l<h^{-\frac{2\alpha^*}{s-1/2}}}\eta_{l}^{1-s}\|e(\eta_l)\|_{H^{1/2}(\Gamma)} &\lesssim k\sum_{1\le \eta_l<h^{-\frac{2\alpha^*}{s-1/2}}}\eta_{l}^{-\epsilon/2}h^{2\alpha^*}\|g\|_{H^{-1/2}(\Gamma)} \\
        &\lesssim \frac{1}{\epsilon}h^{2\alpha^*}\|g\|_{H^{-1/2}(\Gamma)}.
    \end{align*}
        and \eqref{e:subcase} is proven in that subcase.
        \item Case $\alpha_{\Gamma}<1/2$ and $\alpha_{\Gamma}\ge s-1/2$.
        In this case, we choose
        $$\delta=3/2-s+\epsilon/2\in[1/2,1) \qquad \text{and} \qquad \tau=1-s+\epsilon/2\in[0,1/2]$$
        so that $\sigma_1=\sigma_2=2s-1-\epsilon$. Whence, from \eqref{eqn:w_w_h_mu12} we get
        \begin{equation*} \label{eqn:err_w_case2}
        \|e(\eta_l)\|_{H^{1/2}(\Gamma)} \lesssim \eta_l^{s-1 -\epsilon/2}h^{2s-1-\epsilon}\|g\|_{H^{-1/2}(\Gamma)}\le \eta_l^{s-1 -\epsilon/2}h^{2\alpha^*}\|g\|_{H^{-1/2}(\Gamma)}
        \end{equation*}
        and we proceed as in the previous case to obtain \eqref{e:subcase}.
        \item Case $\alpha_{\Gamma}<1/2$ and $\alpha_{\Gamma}< s-1/2$. 
        This time we take
        $$\delta=1-\alpha_{\Gamma}\in[1/2,1) \qquad \text{and} \qquad \tau=1/2-\alpha_{\Gamma}\in[0,1/2].$$
        This gives $\sigma_1=\sigma_2=2\alpha_{\Gamma}$ and thus \eqref{eqn:w_w_h_mu12} becomes
        \begin{equation*} \label{eqn:err_w_case3}
        \|e(\eta_l)\|_{H^{1/2}(\Gamma)} \lesssim \eta_l^{\alpha_{\Gamma}-1/2}h^{2\alpha_{\Gamma}}\|g\|_{H^{-1/2}(\Gamma)} \le \eta_l^{s-1-\widetilde\epsilon/2}h^{2\alpha_{\Gamma}}\|g\|_{H^{-1/2}(\Gamma)}, 
        \end{equation*}
        where $\widetilde\epsilon:=2s-1-2\alpha_{\Gamma}>0$ because $\alpha_{\Gamma}-1/2<s-1$. This implies that
        \begin{align*}
        \frac{\sin(\pi s)}{\pi}k\sum_{1\le \eta_l<h^{-\frac{2\alpha^*}{s-1/2}}}\eta_{l}^{1-s}\|e(\eta_l)\|_{H^{1/2}(\Gamma)} &\lesssim \widetilde\epsilon^{-1}h^{2\alpha_{\Gamma}}\|g\|_{H^{-1/2}(\Gamma)} \\
        &\lesssim \frac{1}{\epsilon}h^{2\alpha^*}\|g\|_{H^{-1/2}(\Gamma)},
    \end{align*}
    which is \eqref{e:subcase}.
    \end{enumerate}
    Gathering the above estimates for all three cases we have
    \begin{equation}
        \label{e:second}
    \|\pi_h \widetilde \psi_k-\dot \psi_{k,h}\|_{H^{1/2}(\Gamma)} \lesssim \frac 1 \epsilon h^{2\alpha^*} \|g\|_{H^{-1/2}(\Gamma)} = \frac{1}{\epsilon}h^{\widetilde \gamma-\epsilon}\|g\|_{H^{-1/2}(\Gamma)}.
    \end{equation}
\end{itemize}
    Estimates \eqref{e:first} and \eqref{e:second} in \eqref{e:des} yield the desired estimate.
\end{proof}

We now extend the stability provided by Lemma 4.1 in \cite{BL2022} to less regular right-hand side. The proof follows that in \cite{BL2022}.

\begin{lemma} \label{prop_stability_Qks}
    Suppose that $g\in H^{-r}(\Gamma)$ and let $v_k:=\mathcal{Q}_k^{-s}(\mathcal{L})g$. Then for any $\beta<2s-r$ we have
    $$\|v_k\|_{\dot H^{\beta}(\Gamma)}\lesssim \max\left\{\frac{1}{1-s},\frac{1}{2s-r-\beta}\right\}\|g\|_{\dot H^{-r}(\Gamma)}.$$
    Consequently, there holds
    $$\|v_k\|_{H^{\beta}(\Gamma)}\lesssim \max\left\{\frac{1}{1-s},\frac{1}{2s-r-\beta}\right\}\|g\|_{H^{-r}(\Gamma)}$$
    provided $\beta,r \in [-1-\alpha_\Gamma,1+\alpha_\Gamma]$.
\end{lemma}

\begin{proof}
    From the definition of $v_k$ (recall that $\eta=e^y$), we deduce that
    \begin{equation}\label{e:stab_Q_int}
    \|v_k\|_{\dot H^{\beta}(\Gamma)}\le \frac{k\sin(\pi s)}{\pi}\sum_{l=-\texttt{M}}^{\texttt{N}}e^{(1-s)y_l}\|w(y_l)\|_{\dot H^{\beta}(\Gamma)}, \quad w(y_l)=(e^{y_l}I+\mathcal{L})^{-1}g.
    \end{equation}
    We now estimate the $\dot H^{\beta}(\Gamma)$ norm of $w(y)$ for any $y\in\mathbb{R}$. Recall that $\{(\eig_j,e_j)\}_{j=1}^{\infty}$ denotes the eigenpairs of $\mathcal{L}$ and set $v:=\mathcal{L}^{-r/2}g\in L^2(\Gamma)$ so that
    \begin{align}
        \|w(y)\|_{\dot H^{\beta}(\Gamma)}^2 &= \|\mathcal{L}^{\frac{\beta}{2}}(e^{y_l}I+\mathcal{L})^{-1}\mathcal{L}^{\frac{r}{2}}v\|_{L^2(\Gamma)}^2 \nonumber \\
        &= \sum_{j=1}^{\infty}\eig_j^{\beta}(e^y+\eig_j)^{-2}\eig_j^r|v_j|^2= \sum_{j=1}^{\infty}\left(\frac{\eig_j^{\frac{\beta+r}{2}}}{e^y+\eig_j}\right)^2|v_j|^2, \label{eqn:inner_w_y}
    \end{align}
    where $v_j:=\langle v,e_j\rangle_{L^2(\Gamma)}$ for $j\ge 1$. Since  $\beta+r<2s<2$ and $\eig_j\ge 1$  for any $j\ge 1$, we have 
    $$\frac{\eig_j^{\frac{\beta+r}{2}}}{e^y+\eig_j}\le \eig_j^{\frac{\beta+r}{2}-1}\le 1.$$
    Alternatively, Young's inequality $ab\le a^p/p+b^q/q$ with $p=1/(1-\theta)$, $q=1/\theta$, and $\theta=(\beta+r)/2$, ensures that
    \begin{equation} \label{eqn:Young}
       \frac{\eig_j^{\frac{\beta+r}{2}}}{e^y+\eig_j}=e^{\left(\frac{\beta+r}{2}-1\right)y}\frac{e^{\left(1-\frac{\beta+r}{2}\right)y}\eig_j^{\frac{\beta+r}{2}}}{e^y+\eig_j}\le e^{\left(\frac{\beta+r}{2}-1\right)y}.
    \end{equation}
Whence, returning to \eqref{e:stab_Q_int} we obtain
    $$\|v_k\|_{\dot H^{\beta}(\Gamma)}\le \frac{k\sin(\pi s)}{\pi}\left[\sum_{l=-\texttt{M}}^0e^{(1-s)y_l}+\sum_{l=1}^{\texttt{N}}e^{\left(\frac{\beta+r}{2}-s\right)y_l}\right]\|v\|_{L^2(\Gamma)},$$
    which is the first desired estimate in disguised upon noting that $\|v\|_{L^2(\Gamma)}=\|g\|_{\dot H^{-r}(\Gamma)}$.

    The second estimate directly follows from the equivalence~\eqref{e:equiv_Hs}.
\end{proof}

\bibliographystyle{plain}
\bibliography{bibliography}

\end{document}